\documentclass[11pt,reqno]{amsart}
\usepackage{amssymb,amsmath,amsthm,mathtools,mathrsfs,slashed}
\usepackage[a4paper, margin=1in, includeheadfoot]{geometry}
\usepackage[inline]{enumitem}
\usepackage{hyperref}
\usepackage[hyperpageref]{backref}
\usepackage[msc-links, abbrev, backrefs]{amsrefs} 
\BibSpecAlias{misc}{thesis}

\usepackage{tikz}
\usetikzlibrary{cd, arrows}
\usepackage[capitalise]{cleveref} 

\usepackage[draft,todonotes={textsize=scriptsize}]{changes}
\definechangesauthor[name=Ai, color=blue]{1}
\definechangesauthor[name=Jiang, color=orange]{2}
\definechangesauthor[name=Jost, color=red]{3}
\setauthormarkuptext{name}

\newcommand{\ad}{\mathrm{ad}\,}
\newcommand{\A}{{A_{\mathrm{ref}}}}

\newcounter{stepnum}
\newcommand{\step}{%
  \par
  \refstepcounter{stepnum}%
  \noindent\textsc{Step \arabic{stepnum}}.\enspace\ignorespaces
}
\crefname{stepnum}{\textsc{Step}}{\textsc{Steps}}
\Crefname{stepnum}{\textsc{Step}}{\textsc{Steps}}

\numberwithin{table}{section}
\numberwithin{equation}{section}

\theoremstyle{plain}
\newtheorem{thm}{Theorem}[section]
\crefname{thm}{Thm.}{Thms.}
\Crefname{thm}{Thm.}{Thms.}
\newtheorem*{thm*}{Theorem}
\newtheorem{lem}[thm]{Lemma}
\crefname{lem}{Lem.}{Lems.}
\Crefname{lem}{Lem.}{Lems.}
\newtheorem{prop}[thm]{Proposition}
\crefname{prop}{Prop.}{Props.}
\Crefname{prop}{Prop.}{Props.}
\newtheorem{cor}[thm]{Corollary}
\crefname{cor}{Cor.}{Cors.}
\Crefname{cor}{Cor.}{Cors.}
\newtheorem*{clm}{Claim}
\theoremstyle{definition}

\crefname{defn}{Defn.}{Defns.}
\Crefname{defn}{Defn.}{Defns.}
\newtheorem{exmp}[thm]{Example} 
\crefname{exmp}{Ex.}{Exs.}
\Crefname{exmp}{Ex.}{Exs.}
\newtheorem{rmk}[thm]{Remark}
\crefname{rmk}{Rem.}{Rem.}
\Crefname{rmk}{Rem.}{Rem.}
\crefname{section}{Sect.}{Sects.}
\Crefname{section}{Sect.}{Sects.}

\usepackage{color}

\title{Variational aspects of the generalized Seiberg--Witten functional}
\date{\today}
\author{Wanjun Ai}
\email{wanjunai@swu.edu.cn}
\address{School of Mathematics and Statistics\\
  Southwest University\\
  Chongqing 400715\\
  P. R. China
}
\address{Max Planck Institute for Mathematics in the Sciences\\
  Inselstr. 22\\
  04103, Leipzig\\
Germany}
\author{Shuhan Jiang}
\email{Shuhan.Jiang@mis.mpg.de}
\address{Max Planck Institute for Mathematics in the Sciences\\
  Inselstr. 22\\
  04103, Leipzig\\
Germany}
\author{Jürgen Jost}
\email{jjost@mis.mpg.de}
\address{Max Planck Institute for Mathematics in the Sciences\\
  Inselstr. 22\\
  04103, Leipzig\\
Germany}
\subjclass[2020]{58E15
}
\keywords{Seiberg--Witten functional, Kapustin--Witten functional, Regularity of weak solutions, Palais--Smale compactness}
\begin{document}
\allowdisplaybreaks
\maketitle
\begin{abstract}
  In this paper, as a step towards a unified mathematical treatment of the
  gauge functionals from quantum field theory that have found profound
  applications in mathematics, we generalize the Seiberg--Witten functional
  that in particular includes the Kapustin--Witten functional as a special
  case. We first demonstrate the smoothness of weak solutions to this generalized
  functional. We then establish the existence of weak solutions under the assumption that the structure group of the bundle is abelian, by verifying the Palais--Smale compactness.
\end{abstract}

\section{Introduction}\label{sec:intro}

Gauge theory has proven to be one of the most powerful tools in the study of low dimensional smooth manifolds, ever since Donaldson's seminal work on anti-self duality equations \cite{Donaldson;Gauge;;1983}. Later, based on a better understanding of the infrared behavior of the $N = 2$ supersymmetric gauge theory  in the nineties \citelist{\cite{Seiberg--Witten;Electromagnetic;Phys;1994}\cite{Seiberg--Witten;QCD;Phys;1994}}, Witten was able to introduce a new set of equations with a simpler gauge group $\mathrm{U}(1)$, which offered a shortcut to many important results of the Donaldson theory \cite{Witten;Monopoles;Math;1994}. These equations are known as the Seiberg--Witten equations. They are first-order equations on a closed oriented Riemannian $4$-manifold $M$:
\begin{align}\label{SWeqs}
  F_A^+-  \frac{i}{2}  \langle e_{\mu}e_{\nu} \sigma, \sigma \rangle e^{\mu} \wedge e^{\nu} = 0, \quad
  \slashed{D}_A \sigma = 0,
\end{align}
where $A$ is a connection of the determinant line bundle of a $\mathrm{spin}^c$-structure $P$ over $M$, $\sigma$ is a section of the half spinor bundle $S^+$ of $P$, and $\slashed{D}_A$ is the twisted Dirac operator of $S^+$.  

The Kapustin--Witten equations are another set of first-order equations on $M$, parameterized by $t \in \mathbb{CP}^1$:
\begin{align}\label{KWeqs}
  (F_A-  \frac{1}{2}  [\sigma, \sigma] + t d_A \sigma)^+= 0, \quad
  (F_A-  \frac{1}{2}  [\sigma, \sigma] - t^{-1} d_A \sigma)^-= 0, \quad
  d_A^* \sigma =0, 
\end{align}
where $A$ is a connection of a principal $\mathrm{SU}(2)$-bundle $P$ over $M$ and  $\sigma$ is an $\ad P$-valued $1$-form on $M$. They come from a topological twisting of the $N=4$ super Yang--Mills theory and are intimately connected to the geometric Langlands program and non-abelian Hodge theory \citelist{\cite{Kapustin;geometric-Langlands-program;;2007}\cite{Liu-Rayan-Tanaka;Kapustin--Witten;;2022}}. For $t=\pm i$, the Kapustin--Witten equations take the following simpler form
\begin{align}\label{KWeqs_t=i}
  F_A - \frac{1}{2}[\sigma,\sigma]=0,\quad (d_A + d_A^*) \sigma =0.
\end{align}
Resembling each other, \eqref{SWeqs} and \eqref{KWeqs_t=i} are special cases of the following first order equations on a closed Riemannian $n$-manifold $M$ equipped with a so-called $\mathrm{spin}^G$-structure $P$:
\begin{align}\label{GSWeqs}
  F_A - \frac{1}{2}\tau(\sigma)=0,\quad \slashed{D}_A \sigma =0,
\end{align}
where $A$ is a connection of a principal $(G/\mathbb{Z}_2)$-bundle $P_{G/\mathbb{Z}_2}$ over $M$ determined by the $\mathrm{spin}^G$-structure $P$, $\sigma$ is a section of an associated vector bundle $S$ of $P$, $\slashed{D}_A$ is a generalized Dirac operator on $S$, and $\tau$ is a quadratic map from $\Gamma(S)$ to $\Omega^2(\mathrm{ad}P_{G/\mathbb{Z}_2})\mathpunct{:}=\Gamma(\Lambda^2 T^*M \otimes \ad P_{G/\mathbb{Z}_2})$.\footnote{One can also replace $F_A$ with its self-dual part $F_A^+$ in \eqref{GSWeqs} when $n=4$.} \eqref{GSWeqs} are known as the generalized Seiberg--Witten equations, which were introduced and studied in \citelist{\cite{Taubes;Nonlinear-3-Dirac;;1999}\cite{Pidstrigach;Seiberg--Witten;;2004}\cite{Haydys;GSW-hyperKähler;;2006}\cite{Callies;GSW-CSD;;2010}\cite{Schumacher;GSW-Swann;;2010}}.

Just like the Seiberg--Witten equations \cite{Jost-Peng-Wang;variational;calc;1996}, \eqref{GSWeqs} satisfy more general second-order equations. To see this, we simply square the LHS of \eqref{GSWeqs} to obtain the following functional
\begin{equation}\label{eq:GSW}
  \mathcal{GSW}(A,\sigma) =  \int_M d\mathrm{vol} \left(|F_A- \frac{1}{2}\tau(\sigma)|^2 + |\slashed{D}_A \sigma|^2  \right). 
\end{equation}
The solutions to \eqref{GSWeqs} constitute the global minimizers of $\mathcal{GSW}(A,\sigma)$. Therefore, they automatically solve the Euler--Lagrange equations of $\mathcal{GSW}(A,\sigma)$
\begin{align}  
  d_A^* F_A + \mathfrak{S}_A(\sigma)&=0, \label{eq:A} \\
  \nabla_A^*\nabla_A \sigma +  \mathfrak{R}_M(\sigma) + \frac{1}{2} \mathfrak{T}(\sigma)&=0,  \label{eq:sigma}
\end{align}
where $\mathfrak{S}_A(\sigma)$ is an $\ad P_{G/\mathbb{Z}_2}$-valued $1$-form of second order in $\sigma$,  $\mathfrak{R}_M$ is an endomorphism of $S$ dependent on the Riemannian curvature of $M$,
and $\mathfrak{T}(\sigma) \in \Gamma(S)$ is of third order in $\sigma$. The
variational structures of the second-order equations \eqref{eq:A} and
\eqref{eq:sigma} provide an alternative strategy, based on methods of
geometric analysis, for constructing solutions of the first-order equations \eqref{GSWeqs}, compared to the index theorems or methods from algebraic geometry. 

The main results of our paper are two-fold: 

First, we prove the smoothness of weak solutions to \eqref{eq:A} and \eqref{eq:sigma}. This is important because for the existence of solutions, one usually needs to turn to the second order equations which naturally arise as the Euler--Lagrange equations of the functional in question. Variational methods then give weak solutions, and one needs to show that these are smooth to produce ordinary solutions. 
\begin{thm}\label{thm:regularity}
  Suppose $(A,\sigma)\in W^{1,2}(\Omega^1(\ad P_{G/\mathbb{Z}_2}))\times W^{1,2}(\Gamma(S))$ is a weak solution of \eqref{eq:A} and \eqref{eq:sigma}, then there exists a gauge transformation $g$, such that  $g(A,\sigma)\mathpunct{:}=(g(A),g^{-1}\sigma)$ is smooth.
\end{thm}
The weak solution $(A,\sigma)$ of $\mathcal{GSW}$ can be described as a coupled system of partial differential equations. In the case where $A$ is a flat connection, the equation for $\sigma$ can be understood as a generalization of harmonic maps, which is evident from the functional definition. On the other hand, when $\sigma$ vanishes, the equation for $A$ reduces to the Yang--Mills equation. Similar to other coupled systems such as the Dirac-harmonic maps \cite{Chen-Jost-Li-Wang;Regularity}, the Yang--Mills--Higgs fields \cite{Ai-Song-Zhu;Yang--Mills--Higgs} and Yang--Mills--Higgs--Dirac fields \cite{Jost-Kessler-Wu-Zhu;Geometric}, this coupled system is critical in the sense that the routine bootstrap argument cannot improve its regularity. 

Since the absolute minimizer of $\mathcal{GSW}$ automatically solves the corresponding first-order equation \eqref{GSWeqs}, our theorem generalizes the regularity results of first order Kapustin--Witten equations \cite{Gagliardo-Uhlenbeck;Kapustin--Witten;Fixed;2012}. Moreover, this full weak regularity theorem allows us to recover classical removal of isolated singularity theorems for Seiberg--Witten equations \cite{Salamon;Removable;;1996} and Yang--Mills equations \cite{Uhlenbeck;Removable;Comm}. This result has important implications for the regularity of solutions to various equations in mathematical physics.

When $\mathcal{GSW}$ is the classical Seiberg--Witten functional, the regularity of weak solutions is proved in \cite{Jost-Peng-Wang;variational;calc;1996}. The proof of our regularity theorem relies on several key ingredients. First, the key estimate in \cite{Jost-Peng-Wang;variational;calc;1996}*{Lem.~3.2} regarding the $L^\infty$ norm of the field $\sigma$, which can be traced back to Taubes \cite{Taubes;equivalence;Math;1980}, does not hold for our more general functional $\mathcal{GSW}$ due to the loss of positivity in the fourth-order term. To overcome this issue, we replace the argument with a classical Moser iteration, which simplifies the proof and overcomes the loss of positivity. Then, the $L^\infty$ boundedness of $\sigma$ results in a growth estimate for connections, similar to the regularity for weak solutions of Yang--Mills fields \cite{Riviere;variations-Yang--Mills-Lagrangian;;2020}. Finally, the Adams-Sobolev embeddings \cite{Adams;Riesz-potentials;Duke;1975} are employed to improve the regularity of the connection beyond the critical index.

The second-order equations might have more solutions than the first order
ones. Such additional solutions would then not be minimizers, and
consequently, saddle point arguments would be needed to show their existence. 
As a first step towards proving such existence results, we  show that the functional $\mathcal{GSW}(A,\sigma)$ satisfies the Palais--Smale condition. This is done under the assumptions that $G$ is abelian and $|\sigma|^2\leq\lambda_0\lvert \tau(\sigma) \rvert$ for some constant $\lambda_0>0$. To make our proof more accessible, we focus on a specific functional, $\mathcal{H}$, defined in \cref{hwf}
However, our existence theorem applies to any functional that satisfies the aforementioned assumptions. In our analysis of the $\mathcal{H}$ functional, we will consider the space $\mathscr{A}^{1,2}\mathpunct{:}=\left\{ \A+a \mathpunct{:} a\in W^{1,2}(\Omega^1(\ad P_{G/\mathbb{Z}_2})) \right\}$, as well as the closed sub-manifold $\mathcal{V}$, defined as 
\[
  \mathcal{V}=\Bigl\{ 
    (A,\sigma)\in \mathscr{A}^{1,2}\times W^{1,2}(\Omega^1(E)) \mathpunct{:}
      \lvert \sigma \rvert^2(x)\leq \lambda_0 \lvert \tau(\sigma) \rvert(x)
      \text{ for all }x\in M \Bigr\},
\]
where $E$ is a complex line bundle, and $\lambda_0>0$ is a constant.
\begin{thm}\label{thm:main}
  The functional $\mathcal{H}$ satisfies the following Palais--Smale condition on $\mathcal{V}$: for any sequence $(A_n,\sigma_n)\in \mathcal{V}$ satisfying 
  \begin{enumerate}
    \item $d \mathcal{H}(A_n,\sigma_n)\to0$ strongly in $W^{-1,2}(\Omega^1(\ad P_{G/\mathbb{Z}_2}))\times W^{-1,2}(\Omega^1(E))$;
    \item $\mathcal{H}(A_n,\sigma_n)\leq c$, for $n=1,2,\ldots$,
  \end{enumerate}
  there exists a subsequence, for simplicity still represented by the index $n$, and a sequence of gauge transformation $g_n\in \mathscr{G}^{2,2}$, such that 
  \[
    g_n(A_n,\sigma_n)\mathpunct{:}=(g_n(A_n),g_n^{-1}\sigma_n)\to (A_\infty,\sigma_\infty)
  \]
  strongly in $W^{1,2}(\Omega^1(\ad P_{G/\mathbb{Z}_2}))\times W^{1,2}(\Omega^1(E))$ as $n\to\infty$; Moreover $(A_\infty,\sigma_\infty)$ is a critical point of $\mathcal{H}$ and 
  \[
    \mathcal{H}(A_\infty,\phi_\infty)=\lim_{n\to\infty}\mathcal{H}(g_n(A_n),g_n^{-1}\sigma_n).
  \]
\end{thm}

It is easy to show that $\mathcal{V}$ is weakly closed and $\mathcal{H}$ is sequentially weakly lower semi-continuous on $\mathcal{V}$ with respect to $W^{1,2}(\Omega^1(\ad P_{G/\mathbb{Z}_2}))\times W^{1,2}(\Omega^1(E))$. Thus, as a corollary, we have
\begin{cor}\label{cor:existence}
  The functional $\mathcal{H}$ is bounded from below on $\mathcal{V}$ and attains its infimum in $\mathcal{V}$. In particular, there exists a weak solution of the Euler--Lagrange equations of $\mathcal{H}$, see \eqref{eq:EL-HW}, which is also smooth according to \cref{thm:regularity}.
\end{cor}

The first challenge in verifying the Palais--Smale condition for the generalized Seiberg--Witten functional is the non-smoothness of the action of gauge transformations on the affine connection space $\mathscr{A}^{1,2}$ when the structure group $G$ is non-abelian. Additionally, the non-abelian structure group makes it difficult to establish the coerciveness of the functional, as we need to control the $W^{1,2}$ norm of connections by the $L^2$ norm of curvature in the functional, see the proof of \cref{lem:coercive}. This becomes problematic when there is no global Coulomb gauge. To overcome this, we restrict ourselves to functionals with an abelian structure group.

Even for the abelian structure group case, we do not have the coerciveness of $\mathcal{GSW}$ on the entire admissible space of weak solutions, due to the absence of certain positive terms that depend on the fourth power of the section $\sigma$, in contrast to the classical Seiberg--Witten case \cite{Jost-Peng-Wang;variational;calc;1996}. To address this, we introduce the closed sub-manifold $\mathcal{V}$ of the Banach space of weak solutions, which preserves the positivity and coerciveness of $\mathcal{H}$ while also being weakly closed and preserving the weak lower semi-continuity of $\mathcal{H}$. By using classical variational methods, we can then establish the existence of weak solutions for $\mathcal{H}$, as stated in \cref{cor:existence}.

We would like to conclude the introduction by some remarks relating to the condition defining $\mathcal{V}$.
\begin{rmk}
  \begin{enumerate}
    \item The condition defining $\mathcal{V}$ simply asserts that
      $\tau(\sigma)(x)\neq0$ whenever $\sigma(x)\neq0$, due to the homogeneity of $\tau(\sigma)$. This condition is essential in providing positive terms that can be used to control the $W^{1,2}$ norms of $\sigma$, and it holds automatically for the Seiberg--Witten functional $\mathcal{SW}$. However, it is not valid for the $\mathcal{KW}$ functional for $G=\mathrm{SU}(2)$;
    \item Moreover, the condition defining $\mathcal{V}$ is also essential in proving the strong convergence of connections, see the \cref{step:clm,step:An-convergence} in the proof of \cref{thm:main};
    \item It will be clear from \cref{rmk:positive-curvature-implies-zero-sigma} that we cannot expect the curvature term in the equivalent form of $\mathcal{GSW}$, see \eqref{gsw2}, to be positive. In that case, we would propose a positive curvature condition that would force the section $\sigma$ to vanish.
  \end{enumerate}
\end{rmk}
\begin{rmk}
The condition for defining $\mathcal{V}$ can be applied to a general functional $\mathcal{GSW}$, and it will become evident from our proof of \cref{thm:main} that the results hold for any functional $\mathcal{GSW}$ as long as the coerciveness condition is met. This condition, in turn, depends on the availability of a global Coulomb gauge. We intend to investigate the existence of this gauge for a non-abelian structure group in a future paper.
\end{rmk}

The rest of the paper will be organized as follows: In \cref{sec:setting}, we provide the geometric settings for the generalized Seiberg--Witten functional, which is modeled on a Dirac $G$-bundle. In \cref{sec: GSW}, we will transform the $\mathcal{GSW}$ into an equivalent form via the Weitzenb\"ock formula, which is more convenient to work with in the calculus of variations. In \cref{sec:regularity}, we will prove the regularity \cref{thm:regularity}. Finally, in \cref{sec:compactness}, we show that the Palais--Smale condition holds for the functional $\mathcal{H}$, and finish the proof of \cref{thm:main} and \cref{cor:existence}.

\section{The geometric settings}\label{sec:setting}

In this section, we provide the geometric setup of the generalized Seiberg--Witten functional. The key notion is that of a Dirac $G$-operator $\slashed{D}_A$, which generalizes the usual twisted Dirac operator and the Dirac operator $d_A + d_A^*$ on Lie algebra valued differential forms. To define $\slashed{D}_A$, we introduce the concepts of Dirac $G$-bundles and $\mathrm{spin}^G$-structures, which can be viewed as the generalizations of the usual Dirac bundles and spin structures. We end this section by proving a Weitzenb\"ock formula for $\slashed{D}_A$. This formula is crucial for the construction of an equivalent form of the generalized Seiberg--Witten functional \eqref{eq:GSW}.
\subsection{Dirac G-bundles}

Let $\mathrm{Cl}(n)$ denote the Clifford algebra of $\mathbb{R}^n$ and $\mathrm{End}(m)$ denote the endomorphism algebra of $\mathbb{R}^m$, i.e., the algebra of $m\times m$ matrices. Let $G \subset \mathrm{End}(m)$ be a (compact) matrix group containing $-1$. There exists a canonical $\mathrm{SO}(n) \times (G/\mathbb{Z}_2)$-action $\rho_{n,G}$ on $\mathrm{Cl}(n)\otimes \mathrm{End}(m)$, which is induced by the canonical action of $\mathrm{SO}(n)$ on $\mathbb{R}^n$ and the adjoint action of $G$ on $\mathrm{End}(m)$.
Let $M$ be an (oriented) Riemannian manifold of dimension $n$. Let $E$ be a vector bundle over $M$ of rank $m$ with structure group $G/\mathbb{Z}_2$. The frame bundle of $TM\otimes E$ admits a $\mathrm{SO}(n) \times (G/\mathbb{Z}_2)$-subbundle $P(E,M) = P_{\mathrm{SO}(n)} \times_M P_{G/\mathbb{Z}_2}$, where $P_{\mathrm{SO}(n)}$ is the frame bundle of $M$ and $P_{G/\mathbb{Z}_2}$ is the principal $G/\mathbb{Z}_2$-bundle encoding the $G/\mathbb{Z}_2$-structure of $E$. We define the Clifford $G$-bundle of the pair $(E,M)$ to be the algebra bundle
\begin{align*}
  \mathrm{Cl}^G(E,M)=P(E,M)\times_{\rho_{n,G}} \left(\mathrm{Cl}(n)\otimes \mathrm{End}(m)\right).
\end{align*}
Note that $\mathrm{Cl}^G(E,M) \cong \Lambda (TM) \otimes \mathrm{End}(E)$ as vector bundles.
The Levi-Civita connection $\nabla$ on $TM$ together with a connection $1$-form $A$ on $P_{G/\mathbb{Z}_2}$ determines a connection on $\mathrm{Cl}^G(E,M)$, which we denote by $\nabla_A$.
Let $S$ be a vector bundle over $M$ equipped with a bundle metric $\langle \cdot, \cdot \rangle$. With a slight abuse of notation, we use $\nabla_A$ once again to denote a metric connection on $S$.
We call $S$ a Dirac $G$-bundle over $M$ if it is a bundle of left modules over $\mathrm{Cl}^G(E,M)$ such that at each $x \in M$,
\begin{align}\label{DG1}
  \langle (e \otimes g) \sigma_1, (e \otimes g) \sigma_2 \rangle = \langle \sigma_1, \sigma_2 \rangle
\end{align}
for all $\sigma_1, \sigma_2 \in S_x$ and all unit vectors $e \in T_xM$ and $g \in G \subset \mathrm{End}(E_x)$, and
\begin{align}\label{DG2}
  \nabla_A(\varphi \sigma) = (\nabla_A \varphi)\sigma + \varphi(\nabla_A \sigma)
\end{align}
for all $\varphi \in \Gamma(\mathrm{Cl}^G(E,M))$ and $\sigma \in \Gamma(S)$.
The Dirac $G$-operator of a Dirac $G$-bundle is a first order differential operator $\slashed{D}_A: \Gamma(S) \rightarrow \Gamma(S)$ defined by 
\begin{align}
  \slashed{D}_A \sigma \mathpunct{:}= \sum_{\mu=1}^n (e_{\mu} \otimes 1) \nabla_{A,e_{\mu}} \sigma
\end{align}
at $x \in M$, where $\{e_{\mu}\}_{\mu=1}^n$ is an orthonormal basis of $T_xM$.
\begin{prop}\label{diracselfadj}
  The Dirac $G$-operator $\slashed{D}_A$ is formally self-adjoint.
\end{prop}
\begin{proof}
  Using \eqref{DG1} and \eqref{DG2}, one can show that
  \begin{align*}
    \langle \slashed{D}_A \sigma_1, \sigma_2 \rangle = \mathrm{div}(V) + \langle \sigma_1, \slashed{D}_A \sigma_2 \rangle,
  \end{align*}
  where $V = - \sum_{\mu=1}^n \langle \sigma_1, (e_{\mu} \otimes 1)\sigma_2 \rangle e_{\mu}$ is a vector field over $M$.
\end{proof}
\begin{exmp}
  The Clifford $G$-bundle $\mathrm{Cl}^G(E,M)$ (viewed as a vector bundle) together with $\nabla_A$ and the bundle metric induced by the standard inner product $\langle \cdot,\cdot \rangle$ on $\mathbb{R}^n$ and the inner product $\langle \cdot,\cdot \rangle$ on $\mathrm{End}(m)$ defined by $\langle A, B \rangle = -\mathrm{Tr}(AB)$ for $A, B \in \mathrm{End}(m)$ is a Dirac $G$-bundle.
\end{exmp}

\subsection{Spin G-structures}

Let $\mathrm{Spin}(n)$ denote the spin group of $\mathbb{R}^n$.
The spin $G$-group $\mathrm{Spin}^G(n)$ is defined as the subgroup of the unit group of $\mathrm{Cl}(n)\otimes \mathrm{End}(m)$ generated by $\mathrm{Spin}(n) \subset \mathrm{Cl}(n)$ and $G \subset \mathrm{End}(m)$.
\begin{lem}\label{spinG}
  The spin $G$-group has the following isomorphism: $\mathrm{Spin}^G(n) \cong \mathrm{Spin}(n) \times_{\mathbb{Z}_2} G \mathpunct{:}= \mathrm{Spin}(n) \times G / \{(1,1),(-1,-1)\}$.
\end{lem}
\begin{proof}
  Consider the map 
  \begin{align*}
    \mathrm{Spin}(n) \times G &\rightarrow \mathrm{Spin}^G(n) \subset \mathrm{Cl}(n) \otimes \mathrm{End}(m)\\
    (a, g) &\mapsto a \otimes g,
  \end{align*}
  which is surjective and has kernel $\{(1,1),(-1,-1)\}$.
\end{proof}
For $G=\mathbb{Z}_2$, $\mathrm{Spin}^{\mathbb{Z}_2}(n) \cong \mathrm{Spin}(n)$. For a general $G$, we have the following commutative diagram
\[ \begin{tikzcd}
  \mathrm{Spin}(n) \times G \arrow{d} \arrow{r}{\mathrm{Ad} \times \mathrm{Id}} & \mathrm{SO}(n) \times G \arrow{d} \\
  \mathrm{Spin}^G(n) \arrow{r}{\mathrm{Ad}^G}& \mathrm{SO}(n) \times (G/\mathbb{Z}_2)
\end{tikzcd} \]   
where the vertical arrows are the canonical quotient maps and the bottom arrow $\mathrm{Ad}^G$ is defined by sending $[a,g] \in \mathrm{Spin}^G(n)$ to $(\mathrm{Ad}(a), [g]) \in \mathrm{SO}(n) \times G/\mathbb{Z}_2$ \footnote{$[a,g]$ is the equivalence class of $(a,g) \in \mathrm{Spin}(n) \times G$ with respect to the quotient in \cref{spinG}, and $[g]$ is the equivalence class of $g \in G$ with respect to the $\mathbb{Z}_2$-quotient of $G$.}.
Therefore there exists a short exact sequence
\begin{align}\label{spinGses}
  \mathrm{Id} \rightarrow \mathbb{Z}_2 \rightarrow \mathrm{Spin}^G(n) \xrightarrow{\mathrm{Ad}^G} \mathrm{SO}(n)\times (G/\mathbb{Z}_2) \rightarrow \mathrm{Id}.
\end{align}

Let $\mathfrak{g}$ denote the Lie algebra of $G$ and $\mathfrak{spin}^G(n)$ denote the Lie algebra of $\mathrm{Spin}^G(n)$. Recall that $\mathfrak{spin}(n)$ can be identified with $\mathrm{Cl}^2(n) \subset \mathrm{Cl}(n)$. 
Under the canonical (vector space) isomorphism $ \mathrm{Cl}(n) \cong \Lambda \mathbb{R}^n $ defined by sending $e_1 \cdots e_r $ to $e_1 \wedge \cdots \wedge e_r$, $\mathfrak{spin}(n) \cong \Lambda^2 \mathbb{R}^n$. 
On the other hand, recall that $\mathfrak{so}(n)$ is the space of skew-symmetric $n \times n$ matrices, which can be also identified with $\Lambda^2 \mathbb{R}^n$ by sending the elementary skew-symmetric $(\mu,\nu)$ matrix to $e_{\mu} \wedge e_{\nu}$. A standard computation shows that the Lie group homomorphism $\mathrm{Ad}: \mathrm{Spin}(n) \rightarrow \mathrm{SO}(n)$ induces a Lie algebra isomorphism
\begin{align}\label{spinso}
  \Lambda^2 \mathbb{R}^n \cong \mathfrak{spin}(n) \rightarrow \mathfrak{so}(n) \cong \Lambda^2 \mathbb{R}^n, \quad  \varphi \mapsto 2 \varphi.
\end{align}

The Lie algebra $\mathrm{Lie}(G/\mathbb{Z}_2)$ of $G/\mathbb{Z}_2$, however, does not have an explicit identification with the Lie algebra $\mathfrak{g}$ of $G$ in general. Let $\pi: G \rightarrow G/\mathbb{Z}_2$ denote the canonical quotient map. We have an abstract isomorphism $d\pi_{\mathrm{Id}}: \mathfrak{g} \rightarrow \mathrm{Lie}(G/\mathbb{Z}_2)$, which is the differential of $\pi$ at the identity.
\begin{lem}\label{liesping}
  The Lie algebra of $\mathrm{Spin}^G(n)$ has the following isomorphism: $\mathfrak{spin}^G(n) \cong \mathfrak{so}(n) \oplus \mathrm{Lie}(G/\mathbb{Z}_2)$, where $\cong $ is defined by sending $(\varphi, \xi)$ to $(2\varphi, d\pi_{\mathrm{Id}}(\xi))$.
\end{lem}
\begin{rmk}
  If $G$ is abelian, one can consider the isomorphism defined by sending $[g] \in G/\mathbb{Z}_2$ to $g^2 \in G$.  $d\pi_{\mathrm{Id}}$ can then be identified with the automorphism of $\mathfrak{g}$ defined by multiplying the factor $2$. We have $\mathfrak{spin}^G(n) \cong \mathfrak{so}(n) \oplus \mathfrak{g}$, where $\cong $ is defined by sending $(\varphi, \xi)$ to $(2\varphi, 2\xi)$.
\end{rmk}

Recall that a spin structure on $M$ is a $\mathbb{Z}_2$-equivariant lift of the frame bundle $P_{\mathrm{SO}(n)}$ to a principal $\mathrm{Spin}(n)$-bundle $P_{\mathrm{Spin}(n)}$ with respect to the double covering $\mathrm{Ad}: \mathrm{Spin}(n) \rightarrow \mathrm{SO}(n)$. Likewise, a $\mathrm{spin}^G$-structure on $M$ is a $\mathbb{Z}_2$-equivariant lift of $P(E,M)$ to a principal $\mathrm{Spin}^G(n)$-bundle $P_{\mathrm{Spin}^G(n)}$ with respect to the double covering $\mathrm{Ad}^G:  \mathrm{Spin}^G(n) \rightarrow \mathrm{SO}(n) \times (G/\mathbb{Z}_2)$. 
Note that  the vector bundle $E$ is part of the data for defining a $\mathrm{spin}^G$-structure on $M$. 
\begin{rmk}
  For $G=\mathrm{U}(1) \subset \mathrm{End}(2)$, a $\mathrm{spin}^G$-structure is known as a $\mathrm{spin}^c$-structure. In this case, $E$ is simply a hermitian line bundle, which is known as the determinant line bundle of the $\mathrm{spin}^c$-structure. 
\end{rmk}   

Let $M$ be a $\mathrm{spin}^G$ manifold. A spinor $G$-bundle $S^G(E,M)$ over $M$ is 
\begin{align*}
  S^G(E,M)=P_{\mathrm{Spin}^G(n)} \times_{\mu_{n,G}} S,
\end{align*}
where $S$ is a left module of $\mathrm{Cl}(n)\otimes \mathrm{End}(m)$ and $\mu_{n,G}$ is the induced left action of $\mathrm{Spin}^G(n)$ on $S$.
A vector $G$-bundle $V^G(E,M)$ over  $M$ is 
\begin{align*}
  V^G(E,M)=P_{\mathrm{Spin}^G(n)} \times_{\mathrm{Ad}} V,
\end{align*}
where $V$ is a bimodule of $\mathrm{Cl}(n)\otimes \mathrm{End}(m)$ and $\mathrm{Ad}$ is the induced adjoint action of $\mathrm{Spin}^G(n)$ on $V$.

\begin{prop}\label{svcliff}
  Both $S^G(E,M)$ and $V^G(E,M)$ are bundles of left modules over $\mathrm{Cl}^G(E,M)$.
\end{prop}
\begin{proof}
  The proof is essentially the same as the proof of \cite{Lawson-Michelsohn;Spin-geometry;;2016}*{Prop.~3.8, Chap.~2}.
\end{proof}

The Levi-Civita connection $1$-form on $P_{\mathrm{SO}(n)}$ together with a connection $1$-form $A$ on $P_{G/\mathbb{Z}_2}$ determines a connection $1$-form on $P_{\mathrm{Spin}^G(n)}$, which again induces a covariant derivative $\nabla_A$ on $S^G(E,M)$ or $V^G(E,M)$.  Supposing that $G$ is compact, one can equip $S$ or $V$ with a $\mathrm{Spin}^G(n)$-invariant inner product to define a bundle metric $\langle \cdot,\cdot \rangle$ on $S^G(E,M)$ or $V^G(E,M)$, which satisfies \eqref{DG1} \cite{Lawson-Michelsohn;Spin-geometry;;2016}. Moreover, $\nabla_{A}$ is a metric connection with respect to $\langle \cdot,\cdot \rangle$ and satisfies \eqref{DG2} by construction. Therefore, both $S^G(E,M)$ and $V^G(E,M)$ can be viewed as Dirac-$G$-bundles, at least when $G$ is compact.

\subsection{Weitzenb\"ock formulas}
In this section, we will derive an abstract Weitzenb\"ock formula that will play a crucial role in the generalized Seiberg--Witten theory. Specifically, we will use it to rewrite the generalized Seiberg--Witten functional into an equivalent form that is suitable for the variational theory.

Let $S$ be a Dirac $G$-bundle. Let $R_S$ denote the curvature $2$-form of the covariant derivative $\nabla_A$ of $S$. 
$R_S$ can be interpreted as a section $\mathfrak{R}_S$ of $\mathrm{End}(S)$ via the formula
\begin{align}
  \mathfrak{R}_S(\sigma) = \frac{1}{2}\sum_{\mu,\nu = 1}^n (e_{\mu} e_{\nu} \otimes R_S(e_{\mu},e_{\nu})) \sigma.
\end{align}
\begin{lem}\label{rsselfadj}
  $\mathfrak{R}_S$ is self-adjoint, i.e., $\langle \sigma_1, \mathfrak{R}_S(\sigma_2) \rangle = \langle  \mathfrak{R}_S(\sigma_1), \sigma_2 \rangle$.
\end{lem}
\begin{proof}
  By \eqref{DG1}, we have
  \begin{align*}
    \langle \sigma_1, (1 \otimes R_S(e_{\mu},e_{\nu}))(\sigma_2) \rangle =  -\langle  (1 \otimes R_S(e_{\mu},e_{\nu}))(\sigma_1), \sigma_2 \rangle
  \end{align*}
  and 
  \begin{align*}
    \langle \sigma_1, (e_{\mu} e_{\nu}  \otimes 1))(\sigma_2) \rangle =  (-1)^2 \langle  (e_{\nu} e_{\mu}  \otimes 1))(\sigma_1), \sigma_2 \rangle = -  \langle  (e_{\mu} e_{\nu}  \otimes 1))(\sigma_1), \sigma_2 \rangle.
  \end{align*}
\end{proof}
Let $\nabla_A^*$ denote the formal adjoint operator of $\nabla_A$.
\begin{prop}[General Weitzenb\"ock formula]\label{gwf}
  $\slashed{D}_A^2 = \nabla_A^* \nabla_A + \mathfrak{R}_S$.
\end{prop}
\begin{proof}
  The proof is essentially the same as the proof of \cite{Lawson-Michelsohn;Spin-geometry;;2016}*{Thm.~8.2, Chap.~2}.
\end{proof}

Let's take $S$ to be $S^G(E,M)$ or $V^G(E,M)$. By \cref{liesping}, one can decompose $R_S$ as
\begin{align}\label{rrf}
  R_S = R_M + F_A,
\end{align}
where $R_M$ is the Riemannian curvature $2$-form and $F_A$ is the curvature $2$-form of $A$. Likewise, we can define
\begin{align}\label{rm}
  \mathfrak{R}_M(\sigma) = \frac{1}{2}\sum_{\mu,\nu = 1}^n (e_{\mu} e_{\nu} \otimes R_M(e_{\mu},e_{\nu})) \sigma,
\end{align}
where we view $R_M(e_{\mu},e_{\nu})$ as an element in $\mathfrak{spin}(n)$ using \eqref{spinso}, and 
\begin{equation}\label{eq:F_A}
  \mathfrak{F}_A(\sigma) = \frac{1}{2}\sum_{\mu,\nu = 1}^n (e_{\mu} e_{\nu} \otimes (d\pi_{\mathrm{Id}})^{-1} (F_A(e_{\mu},e_{\nu}))) \sigma.
\end{equation}
It follows from \cref{gwf} that
\begin{prop}\label{W2}
  For $S=S^G(E,M)$ or $V^G(E,M)$, we have
  \[
    \slashed{D}_A^2 = \nabla_{A}^* \nabla_{A} + \mathfrak{R}_M + \mathfrak{F}_A.
  \]
\end{prop}

We will provide several examples of Weitzenb\"ock formulas in the following, by considering $S$ as various different Dirac $G$-bundles. In these examples, we can write down the abstract terms $\mathfrak{R}_M$ and $\mathfrak{T}_A$ in \cref{W2} in a more concrete form. These will be used to demonstrate that our generalized Seiberg--Witten functional, defined later, encompasses several classical functionals, such as the Seiberg--Witten and Kapustin--Witten functionals.
\begin{exmp}
  Let's take $G=\mathrm{U}(1)$. We can take $S$ to be a (half) spinor bundle associated with the $\mathrm{Spin}^c$-structure.
  \begin{clm}[\cite{Jost;Riemannian-Geometry-geometric-analysis;;2017}*{Thm. 4.4.2}]
    We have 
    \begin{align}\label{ws}
      \slashed{D}_A^2 = \nabla_{A}^* \nabla_{A} + \frac{1}{4} s + \frac{1}{2}F_A,
    \end{align}
    where $s$ is the scalar curvature of $M$, $F_A= \sum_{\mu < \nu} F_A(e_{\mu},e_{\nu}) e_{\mu} \wedge e_{\nu}$ acts on $\sigma$ via the identification $e_{\mu} \wedge e_{\nu} \mapsto e_{\mu}e_{\nu}$.
  \end{clm}
\end{exmp}

\begin{exmp}
  Let's take $S$ to be $V^G(E,M)=\Lambda T^*M \otimes \ad P_{G/\mathbb{Z}_2}$. $V^G(E,M)$ can be identified with the Clifford $G$-bundle $\mathrm{Cl}^G(E,M)$ using the Riemannian metric. 
  \begin{clm}
    We have
    \begin{align}\label{wv}
      \slashed{D}_A^2 = \nabla_{A}^* \nabla_{A} + \mathrm{Ric} + F_A,
    \end{align}
    where $\mathrm{Ric}(\cdot)$ is the so-called Weitzenb\"ock operator and is defined by
    \begin{align}\label{wo}
      \mathrm{Ric}(\sigma)(e_1,\cdots,e_p) = \sum_{\mu=1}^n \sum_{\nu=1}^p (R_M(e_{\mu},e_{\nu})\sigma)(e_1,\cdots,e_{\nu-1}, e_{\mu}, e_{\nu+1}, \cdots,e_p),
    \end{align}
    $F_A$ is defined in a similar manner, i.e.,
    \[
      F_A(\sigma)(e_1,\cdots,e_p) = \sum_{\mu=1}^n \sum_{\nu=1}^p (F_A(e_{\mu},e_{\nu})\sigma)(e_1,\cdots,e_{\nu-1}, e_{\mu}, e_{\nu+1}, \cdots,e_p).
    \]
  \end{clm}
  \begin{rmk}\label{Ric}
    Applying \eqref{wo} to the metric dual $X^{\sharp}=g(X,\cdot)$ of the vector field $X$, we have
    $
    \mathrm{Ric}(X^{\sharp})(Y)=\sum_{\mu=1}^n (R(e_\mu,Y)X^{\sharp})(e_\mu)
    =\sum_{\mu=1}^n g(R(e_\mu,Y)X,e_\mu)
    =-\sum_{\mu=1}^n g(R(e_\mu,Y)e_\mu,X)
    =\mathrm{Ric}(X,Y),
    $
    where $\mathrm{Ric}(\cdot,\cdot)$ is the Ricci curvature tensor of $g$, hence the notation. 
  \end{rmk}
  \begin{proof}
    Using \eqref{DG1} and the identification $\mathrm{Cl}(n) \cong \Lambda \mathbb{R}^n$, we obtain
    \begin{align*}
                        &\mathfrak{R}_M(\sigma)(e_{k_1},\cdots,e_{k_p}) = \frac{1}{2} \sum_{\mu,\nu=1}^n \langle e_{\mu}e_{\nu}(R_M(e_{\mu},e_{\nu})(\sigma)), e^{k_1} \wedge \cdots \wedge e^{k_p}\rangle \\
                        &\qquad=-\frac{1}{2} \sum_{\mu,\nu=1}^n \langle R_M(e_{\mu},e_{\nu})(\sigma), e_{\mu}e_{\nu}(e^{k_1} \wedge \cdots \wedge e^{k_p})\rangle \\
                        &\qquad=-\frac{1}{2}\sum_{\mu, \nu=1}^n\langle R_M(e_{\mu},e_{\nu})(\sigma), e_{\mu}e_{\nu}(e^{k_1}) \wedge \cdots \wedge e^{k_p} + \cdots + e^{k_1} \wedge \cdots \wedge (e_{\mu} e_{\nu})(e^{k_p})\rangle \\
                        &\qquad=\begin{multlined}[t]
                          -\frac{1}{2}\left(\sum_{\nu=1}^n\langle R_M(e_{k_1},e_{\nu})(\sigma), e^{\nu} \wedge \cdots \wedge e^{k_p}  \rangle - \sum_{\mu=1}^n\langle R_M(e_{\mu},e_{k_1})(\sigma), e^{\mu} \wedge \cdots \wedge e^{k_p}  \rangle  + \cdots \right. \\
                          \left. + \sum_{\nu=1}^n\langle R_M(e_{k_p},e_{\nu})(\sigma), e^{k_1} \wedge \cdots \wedge e^{\nu}  \rangle - \sum_{\mu=1}^n\langle R_M(e_{\mu},e_{k_p})(\sigma), e^{k_1} \wedge \cdots \wedge e^{\mu}  \rangle \right)
                        \end{multlined} \\
                        &\qquad= \sum_{\alpha=1}^p \sum_{\mu=1}^n \langle R_M(e_{\mu},e_{k_{\alpha}})(\sigma), e^{k_1} \wedge \cdots \wedge e^{k_{\alpha-1}} \wedge e^{\mu} \wedge e^{k_{\alpha+1}} \wedge  \cdots e^{k_p} \rangle 
                        = \mathrm{Ric}(\sigma).
    \end{align*}
    One can show that $\mathfrak{F}_A(\sigma) = F_A(\sigma)$ using a similar computation as above.
  \end{proof}
\end{exmp}

\begin{exmp}
  When $G$ is $\mathrm{U}(1)$, one can also construct the following Dirac $G$-bundle
  \begin{align*}
    S=H(E,M) \mathpunct{:}= P_{\mathrm{Spin}^c(n)} \times_{\rho_{h}} \mathbb{C}\mathrm{l}(n) \cong \Lambda T^*M \otimes_{\mathbb{R}} E, 
  \end{align*}
  where $\mathbb{C}\mathrm{l}(n)\mathpunct{:}=\mathrm{Cl}(n) \otimes_{\mathbb{R}} \mathbb{C}$ is the complex Clifford algebra and $\rho_{h}$ is defined as 
  \begin{align*}
    \rho_{h}: \mathrm{Spin}^c(n) \times (\mathrm{Cl}(n) \otimes_{\mathbb{R}} \mathbb{C}) &\rightarrow \mathrm{Cl}(n) \otimes_{\mathbb{R}} \mathbb{C} \\
    [a,\exp(is)] \times (\varphi \otimes t) &\mapsto \mathrm{Ad}(a)(\varphi) \otimes \exp(2is)t.
  \end{align*}
It is not difficult to demonstrate that $H(E,M)$ can also be seen as a bundle of left modules over $\mathrm{Cl}^G(E,M)$. This is possible because $G$ is abelian. One can equip $H(E,M)$ with a hermitian bundle metric $\langle \cdot, \cdot \rangle_h$\footnote{Our convention is that $\langle \cdot, \cdot \rangle_h$ is anti-linear in its first argument.} by assigning to its fiber $ \mathbb{C}\mathrm{l}(n)$ the standard hermitian inner product. We also require that $\langle \cdot, \cdot \rangle\mathpunct{:}=\mathrm{Re}(\langle \cdot, \cdot \rangle_h)$ satisfies \eqref{DG1}. $H(E,M)$ is then a Dirac $G$-bundle, which we call of hybrid type. Using a similar argument as in the proof of \eqref{wv}, one can show that the Weitzenb\"ock formula for the Dirac $G$-operator on $H(E,M)$ takes the same form
  $
  \slashed{D}_A^2 = \nabla_{A}^* \nabla_{A} + \mathrm{Ric} + F_A.
  $
\end{exmp}

\section{Generalized Seiberg--Witten functional}\label{sec: GSW}

For a vector bundle $E$ with a fixed metric connection $\A$, we define the
Sobolev space $W^{k,p}(E)$ as the completion of the  space of smooth
sections of $E$ under the norm  defined as 
\[
  \lVert \sigma \rVert_{W^{k,p}}\mathpunct{:}=\left( \sum_{\lvert \alpha \rvert=1}^k \int \lvert \nabla _\A^\alpha \sigma \rvert \right)^{1/p},\quad \sigma\in C^\infty(\Gamma(E)).
\]
The usual Sobolev embedding theorem and H\"older inequalities are valid for $W^{k,p}(E)$. 

In particular, the space of connections we will consider is the affine space of $W^{1,2}$ 1-forms:
\[
  \mathscr{A}^{1,2} \mathpunct{:}= \left\{ \A+a \mathpunct{:}a\in W^{1,2}(\Omega^1(\ad P_{G/\mathbb{Z}_2})) \right\},
\]
and the space of sections of a Dirac $G$-bundle $S$:
\[
  \mathscr{S}^{1,2}\mathpunct{:}= W^{1,2}(S) \mathpunct{:}=\left\{ \sigma\in\Gamma(S)\mathpunct{:} \sigma\in W^{1,2}(S) \right\}.
\]

Let 
$
\tau(\sigma) \mathpunct{:}= \sum_{1 \leq \mu < \nu \leq n} \sum_a \langle (e_{\mu}e_{\nu} \otimes d\pi_{\mathrm{Id}}^{-1}(\xi_a))\sigma, \sigma \rangle (e^{\mu} \wedge e^{\nu}) \otimes \xi_a,
$
where $\{\xi_a\}$ is an orthonormal basis of $\mathrm{Lie}(G/\mathbb{Z}_2)$. The Generalized Seiberg--Witten functional is defined as
\begin{equation}\label{gsw}
  \mathcal{GSW}(A,\sigma) =  \int_M d\mathrm{vol} \left(|F_A- \frac{1}{2}\tau(\sigma)|^2 + |\slashed{D}_A \sigma|^2  \right),\quad (A,\sigma)\in \mathscr{A}^{1,2}\times \mathscr{S}^{1,2}. 
\end{equation}
Now, we will apply the generated Weitzenb\"ock formula to rewrite the functional $\mathcal{GSW}$ in a form that separates the terms of $A$ and $\sigma$, and fits well in the study of the variational properties of $\mathcal{GSW}$.
\begin{prop}
For $(A,\sigma)\in \mathscr{A}^{1,2}\times \mathscr{S}^{1,2}$, the generalized Seiberg--Witten functional $\mathcal{GSW}$ can be written equivalently as:
  \begin{equation}\label{gsw2}
    \mathcal{GSW}(A,\sigma) = \int_M d\mathrm{vol} \left(|\nabla_A \sigma|^2 + |F_A|^2 + \langle \sigma, \mathfrak{R}_M(\sigma)  \rangle  + \frac{1}{4} |\tau(\sigma)|^2 \right).
  \end{equation}
\end{prop}
\begin{proof}
  By \cref{diracselfadj,W2}, we have
  \begin{align*}
    \int_M d\mathrm{vol} \langle \slashed{D}_A \sigma, \slashed{D}_A \sigma \rangle = \int_M d\mathrm{vol} \langle \sigma, \slashed{D}_A^2 \sigma \rangle = \int_M d\mathrm{vol} \left(|\nabla_A \sigma|^2 + \langle \sigma, \mathfrak{R}_M(\sigma)  \rangle  + \langle \sigma, \mathfrak{F}_A(\sigma)  \rangle \right).
  \end{align*}
  The rest of the proof follows from direct computations:
  \begin{align*}
                &\int_M d\mathrm{vol} \left(|F_A- \frac{1}{2}\tau(\sigma)|^2 + |\slashed{D}_A \sigma|^2  \right) = \int_M d\mathrm{vol} \left(|F_A|^2 + \frac{1}{4}|\tau(\sigma)|^2 - \langle F_A, \tau(\sigma) \rangle + |\slashed{D}_A \sigma|^2  \right) \\
                &\qquad= \int_M d\mathrm{vol} \left(|F_A|^2 + \frac{1}{4}|\tau(\sigma)|^2 - \langle \sigma, \mathfrak{F}_A(\sigma)  \rangle+ |\nabla_A \sigma|^2 + \langle \sigma, \mathfrak{R}_M(\sigma)  \rangle  + \langle \sigma, \mathfrak{F}_A(\sigma) \rangle  \right) \\
                &\qquad= \int_M d\mathrm{vol} \left(|\nabla_A \sigma|^2 + |F_A|^2 + \langle \sigma, \mathfrak{R}_M(\sigma)  \rangle  + \frac{1}{4} |\tau(\sigma)|^2 \right),
  \end{align*}
  where we use the definition of $\tau(\sigma)$ and \eqref{eq:F_A} to pass to the second line.
\end{proof}

\begin{prop}
  The Euler--Lagrange equations of the generalized Seiberg--Witten functional $\mathcal{GSW}$ are given by \eqref{eq:A} and \eqref{eq:sigma}, i.e.,
  \begin{align*}  
    d_A^* F_A + \mathfrak{S}_A(\sigma)&=0, \\
    \nabla_A^*\nabla_A \sigma +  \mathfrak{R}_M(\sigma) + \frac{1}{2} \mathfrak{T}(\sigma)&=0,
  \end{align*}
  where
  $
  \mathfrak{T}(\sigma) = \sum_{\mu<\nu} \sum_a  \langle (e_{\mu}e_{\nu} \otimes d\pi_{\mathrm{Id}}^{-1}(\xi_a))\sigma, \sigma \rangle (e_{\mu}e_{\nu} \otimes d\pi_{\mathrm{Id}}^{-1}(\xi_a))\sigma,
  $
  and $\mathfrak{S}_A(\sigma)$ is an $\ad P_{G/\mathbb{Z}_2}$-valued $1$-form defined through the identity $\langle \nabla_A \sigma, B \sigma \rangle = \langle \mathfrak{S}_A(\sigma), B \rangle$, where $B\in \Omega^1(\ad P_{G/\mathbb{Z}_2})$. 
\end{prop}
\begin{rmk}
  Locally, we can write $\langle \nabla_A \sigma, B \sigma \rangle = \sum_{\mu,\nu} \langle e^{\mu} \wedge \nabla_{A,\mu} \sigma, e^{\nu} \wedge B_{\nu} \sigma \rangle  = \sum_{\mu} \langle \nabla_{A,\mu} \sigma, B_{\mu} \sigma \rangle$, where $\nabla_{A,\mu} = \iota_{e_{\mu}} \nabla_A$. $\mathfrak{S}_A(\sigma)$ can be then written as
  $\mathfrak{S}_A(\sigma) = \sum_{\mu} (\sigma \otimes (\nabla_{A, \mu} \sigma)^t)\otimes e^{\mu}$, where $(\nabla_{A, \mu} \sigma)^t\mathpunct{:}= \langle \nabla_{A, \mu} \sigma, \cdot \rangle$.
\end{rmk}
\begin{proof}
  The calculation of the variation of the Yang--Mills term is standard. We have
  \begin{align*}
    \frac{d}{dt}\bigg|_{t=0} \int_M d\mathrm{vol} |F_{A+tB}|^2 = 2 \int_M d\mathrm{vol} \langle F_A, d_A B \rangle =  2 \int_M d\mathrm{vol} \langle d_A^*F_A, B \rangle,
  \end{align*}
  where $d_A^*$ is the formal adjoint operator of $d_A$. For the $|\nabla_A \sigma|^2$-term, we have
  \begin{align*}
    \frac{d}{dt}\bigg|_{t=0} \int_M d\mathrm{vol} |\nabla_{A+tB} \sigma|^2 = 2 \int_M d\mathrm{vol} \langle \nabla_A \sigma, B \sigma \rangle = 2 \int_M d\mathrm{vol} \langle \mathfrak{S}_A(\sigma), B  \rangle, 
  \end{align*}
  and
  \begin{align*}
    \frac{d}{dt}\bigg|_{t=0} \int_M d\mathrm{vol} |\nabla_{A} (\sigma + t \beta)|^2 = 2 \int_M d\mathrm{vol} \langle \nabla_A \sigma, \nabla_A \beta \rangle = 2 \int_M d\mathrm{vol} \langle \nabla_A^* \nabla_A \sigma, \beta  \rangle.
  \end{align*}
  For the other terms, noting that, by \cref{rsselfadj}, we have
  \begin{align*}
    \frac{d}{dt}\bigg|_{t=0}|\tau(\sigma + t \beta)|^2 &= 2 \left\langle \tau(\sigma), \frac{d}{dt}\bigg|_{t=0}\tau(\sigma + t \beta)\right\rangle \\
                                                       &= 4 \langle \tau(\sigma), \sum_{\mu<\nu} \sum_a \langle (e_{\mu}e_{\nu} \otimes d\pi_{\mathrm{Id}}^{-1}(\xi_a))\sigma, \beta \rangle (e^{\mu} \wedge e^{\nu}) \otimes \xi_a\rangle \\
                                                       &=4 \sum_{\mu<\nu} \sum_a   \langle (e_{\mu}e_{\nu} \otimes d\pi_{\mathrm{Id}}^{-1}(\xi_a))\sigma, \beta \rangle  \langle \tau(\sigma),(e^{\mu} \wedge e^{\nu}) \otimes \xi_a\rangle \\
                                                       &=4 \sum_{\mu<\nu} \sum_a  \langle (e_{\mu}e_{\nu} \otimes d\pi_{\mathrm{Id}}^{-1}(\xi_a))\sigma, \beta \rangle  \langle (e_{\mu}e_{\nu} \otimes d\pi_{\mathrm{Id}}^{-1}(\xi_a))\sigma, \sigma \rangle \\
                                                       &= 4 \langle \mathfrak{T}(\sigma), \beta \rangle,
  \end{align*}
  and $\frac{d}{dt}\big|_{t=0} \langle \sigma + t \beta, \mathfrak{R}_M(\sigma+ t \beta)  \rangle = 2 \langle \mathfrak{R}_M(\sigma), \beta\rangle$.
\end{proof}

With this rewritten form of $\mathcal{GSW}$ and its Euler--Lagrange equations, we can easily obtain:
\begin{cor}
  Suppose $\sigma$ is a weak solution to \eqref{eq:sigma}, then 
  \begin{equation}\label{eq:Moser-iteration}
    \frac{1}{2}\Delta\lvert \sigma \rvert^2-\lvert \nabla _A\sigma \rvert^2-\left\langle \mathfrak{R}_M(\sigma),\sigma \right\rangle-\frac{1}{4}\lvert \tau(\sigma) \rvert^2=0,
  \end{equation}
  where $\Delta$ is the (negative) analytic Laplacian.
\end{cor}
\begin{proof}
  This follows by taking the inner product of \eqref{eq:sigma} with $\sigma$, and using
  \[
    \frac{1}{2}\Delta \lvert \sigma \rvert^2=\left\langle \nabla _A\sigma,\nabla _A\sigma \right\rangle-\left\langle \nabla _A^*\nabla _A\sigma,\sigma \right\rangle.
  \]
\end{proof}
\begin{cor}\label{cor:L-infty-boundedness-via-Morse}
  Suppose $B'\subset B\subset M$, where $B'$ and $B$ are balls in $M$ with metric properties close to $B_1^4\subset B_2^4\subset \mathbb{R}^4$. Then, for weak solution of \eqref{eq:sigma}, we have 
  \[
    \sup_{x\in B'}\lvert \sigma(x) \rvert\leq C\lVert \sigma \rVert_{L^2(B)}.
  \]
  In particular, a dilating argument show that this estimates holds on any compact manifold $M$.
\end{cor}
\begin{proof}
  Note that \eqref{eq:Moser-iteration} implies 
  \begin{equation}\label{eq:Bochner}
    \Delta\lvert \sigma \rvert^2\geq 2\left\langle \mathfrak{R}_M(\sigma),\sigma \right\rangle \geq 2\lambda\lvert \sigma \rvert^2,
  \end{equation}
  where $\lambda$ is the lower bound of the Ricci curvature. The result follows from the classical Moser iteration.
\end{proof}
\begin{rmk}\label{rmk:positive-curvature-implies-zero-sigma}
  The \eqref{eq:Bochner} implies that if $M$ has a positive lower   Ricci
  curvature bound, then $\sigma$ vanishes identically on $M$.
\end{rmk}

Next, we will show that the generalized $\mathcal{GSW}$ includes several functionals from quantum field theories, specifically the Seiberg--Witten functional and Kapustin--Witten functional. Since our proof of the existence theorem for the generalized Seiberg--Witten equations requires that the gauge group is abelian, we 
also introduce a functional $\mathcal{H}$ with gauge group $U(1)$. In fact, $\mathcal{H}$ can be viewed as a $\mathrm{SU}(2)$ Yang--Mills functional with an additional gauge fixing term, see \cref{rmk: H}.

\begin{exmp}[Seiberg--Witten functional]


  For simplicity, let's consider a $3$-dimensional \emph{variant} of the original $4$-dimensional Seiberg--Witten functional. Let $(M,g)$ be a $3$-dimensional $\mathrm{spin}^c$  manifold. $\mathrm{Spin}^c(3)$ can be identified as $\mathrm{Sp}(1) \times_{\mathbb{Z}_2} \mathrm{U}(1)$. Let's consider the spinor bundle $S$ over $M$ defined via the representation
  \begin{align*}
    (\mathrm{Sp}(1) \times_{\mathbb{Z}_2} \mathrm{U}(1)) \times \mathbb{H} &\rightarrow \mathbb{H} \\
    (a, b) \times s &\mapsto asb.
  \end{align*}
  To derive an explicit expression for $\tau(s)$, $s = s_1 + s_2 i + s_3 j + s_4 k \in \mathbb{H}$, where $i, j, k$ are the imaginary units of $\mathbb{H}$, we make the following identifications
  \begin{align*}
    e_1 \mapsto i, \quad e_2 \mapsto j, \quad e_3 \mapsto k.
  \end{align*}
  We have
  \begin{align*}
    \tau(s) &= \left\langle e_1e_2 s\frac{i}{2}, s \right\rangle e_1 \wedge e_2 \otimes i + \left\langle e_2e_3 s\frac{i}{2}, s \right\rangle e_2 \wedge e_3 \otimes i + \left\langle e_1e_3 s\frac{i}{2}, s \right\rangle e_1 \wedge e_3 \otimes i \\
            &= \frac{1}{2}\left(\mathrm{Re}(\overline{ksi}s)e_1 \wedge e_2 \otimes i + \mathrm{Re}(\overline{isi}s)e_2 \wedge e_3 \otimes i - \mathrm{Re}(\overline{jsi}s)e_1 \wedge e_3 \otimes i \right)\\
            &= (s_1s_3 -s_2s_4) e_1 \wedge e_2 \otimes i + \frac{1}{2}(-s_1^2 - s_2^2 + s_3^2 + s_4^2) e_2 \wedge e_3 \otimes i + (s_1s_4 + s_2 s_3) e_1 \wedge e_3 \otimes i.
  \end{align*}
  Hence,
  \begin{align*}
    |\tau(s)|^2 &= \frac{1}{4}((s_3^2 + s_4^2) - (s_1^2 +s_2^2))^2 + \left((s_1s_3-s_2s_4)^2 + (s_1s_4 + s_2 s_3)^2\right) \\
                &= \frac{1}{4} ((s_3^2 + s_4^2) - (s_1^2 +s_2^2))^2 + (s_1^2 +s_2^2)(s_3^2 + s_4^2) \\
                &= \frac{1}{4}(s_1^2 +s_2^2 +s_3^2 +s_4^2)^2 
                = \frac{1}{4}|s|^4.
  \end{align*}
  By \eqref{ws}, the $3$-dimensional Seiberg--Witten functional can be written as
  \[
    \mathcal{SW}(A,\sigma) \mathpunct{:}= \int_M d\mathrm{vol} \left(|\nabla_A \sigma|^2 + |F_A|^2 + \frac{s}{4} |\sigma|^2  + \frac{1}{16} |\sigma|^4 \right).
  \]
  For the $4$-dimensional case, one needs to work with the self-dual curvature $2$-form $F_A^+$.
   Similar computations then yield the Seiberg--Witten functional considered in \cite{Jost-Peng-Wang;variational;calc;1996}.
\end{exmp}

\begin{exmp}[Kapustin--Witten functional]
Let $(M, g)$ be an $n$-dimensional Riemannian manifold, typically with $n=4$. Let $P$ be a principal $G$-bundle over $M$. Let's consider the Dirac $G$-bundle $S=\Lambda^* TM \otimes \ad P_{G/\mathbb{Z}_2}$. 
  A standard computation shows that the Dirac $G$-operator on $S$ takes the form
  $
  \slashed{D}_A  = d_A + d_A^*,
  $
  where $d_A: \Omega^p(\ad P_{G/\mathbb{Z}_2}) \rightarrow \Omega^{p+1}(\ad P_{G/\mathbb{Z}_2})$ is the exterior covariant derivative, and $d_A^*$ is the formal adjoint operator of $d_A$.

  Let's derive an explicit expression for $\tau(s)$ for $s \in \Lambda^1 \mathbb{R}^n \otimes \mathrm{Lie}_{G/\mathbb{Z}_2}$. In this case, one can simply write
  \begin{align*}
    \tau(s)= \sum_{\rho,\sigma} \sum_{\mu < \nu} \sum_a \langle (e_{\mu} e_{\nu} \otimes \xi_a) (e^{\rho} \otimes s_{\rho}), e^{\sigma} \otimes s_{\sigma} \rangle e^{\mu} \wedge e^{\nu} \otimes \xi_a,
  \end{align*}
  where $\xi_a \in \mathrm{Lie}_{G/\mathbb{Z}_2}$ acts on $s_{\rho}$ via the adjoint action. We then have
  \begin{align*}
    \tau(s) &=  \sum_{\mu < \nu} \sum_a \left(\langle (e_{\mu} e_{\nu} \otimes \xi_a) (e^{\mu} \otimes s_{\mu}), e^{\nu} \otimes s_{\nu} \rangle + \langle (e_{\mu} e_{\nu} \otimes \xi_a) (e^{\nu} \otimes s_{\nu}), e^{\mu} \otimes s_{\mu} \rangle \right) e^{\mu} \wedge e^{\nu} \otimes \xi_a \\
            &= \sum_{\mu < \nu} \sum_a \left(\langle [\xi_a, s_{\mu}], s_{\nu} \rangle - \langle [\xi_a, s_{\nu}], s_{\mu} \rangle \right) e^{\mu} \wedge e^{\nu} \otimes \xi_a \\
            &= \sum_{\mu < \nu} \sum_a \left(\langle \xi_a, [s_{\mu}, s_{\nu}] \rangle - \langle \xi_a, [s_{\nu}, s_{\mu}] \rangle \right) e^{\mu} \wedge e^{\nu} \otimes \xi_a \\
            &= \sum_{\mu < \nu} e^{\mu} \wedge e^{\nu} \otimes ([s_\mu, s_{\nu}] - [s_\nu, s_{\mu}])=[s,s].
  \end{align*}
  The Kapustin--Witten functional takes the form,
  \[
    \mathcal{KW}(A,\sigma) \mathpunct{:}= \int_M \left(|F_A-[\sigma,\sigma]/2|^2 + |d_A \sigma|^2 + |d_A^* \sigma|^2 \right) d\mathrm{vol},
  \]
  where $\sigma$ is a section of $T^*M \otimes \ad P_{G/\mathbb{Z}_2}$. Note that this is exactly the same as the one considered in \cite{Gagliardo-Uhlenbeck;Kapustin--Witten;Fixed;2012}. Alternatively,
  \[
    \mathcal{KW}(A,\sigma) = \int_M \left(|F_A|^2 + |\nabla_A \sigma|^2 + \langle \sigma, \mathrm{Ric}(\sigma) \rangle + |[\sigma,\sigma]|^2/4 \right) d\mathrm{vol}.
  \]
\end{exmp}
\begin{exmp}[$\mathcal{H}$-functional]\label{hwf}
  Let $(M,g)$ be an $n$-dimensional Riemannian manifold. Let $P$ be a principal $\mathrm{U}(1)$-bundle over $M$. Let $E$ be a hermitian line bundle associated to $P$ via the standard representation of $U(1)$. Let's consider the Dirac $G$-bundle $S=H(E,M)\cong \Lambda T^*M \otimes E$ of hybrid type. The Dirac $G$-operator $\slashed{D}_A$ takes the same form as in the case of Kapustin--Witten theory, i.e., $\slashed{D}_A = d_A + d_A^*$. To write down the action functional, we again need to derive an explicit expression for $\tau(s)$, $s \in \Lambda^1 \mathbb{R}^n \otimes \mathbb{C}$. We have
  \begin{align*}
    \tau(s) &= \sum_{\rho,\sigma} \sum_{\mu < \nu} \langle (e_{\mu} e_{\nu} \otimes i) (e^{\rho} \otimes s_{\rho}), e^{\sigma} \otimes s_{\sigma} \rangle e^{\mu} \wedge e^{\nu} \otimes i \\
            &=  \sum_{\mu < \nu} \left(\langle (e_{\mu} e_{\nu} \otimes i) (e^{\mu} \otimes s_{\mu}), e^{\nu} \otimes s_{\nu} \rangle + \langle (e_{\mu} e_{\nu} \otimes i) (e^{\nu} \otimes s_{\nu}), e^{\mu} \otimes s_{\mu} \rangle \right) e^{\mu} \wedge e^{\nu} \otimes i \\
            &= \sum_{\mu < \nu} \left(\langle i s_{\mu}, s_{\nu} \rangle - \langle i s_{\nu}, s_{\mu} \rangle \right) e^{\mu} \wedge e^{\nu} \otimes i \\
            &= \sum_{\mu < \nu} \mathrm{Re}(\langle is_{\mu}, s_{\nu} \rangle_h - \langle is_{\nu}, s_{\mu} \rangle_h) e^{\mu} \wedge e^{\nu} \otimes i\\
            &= \sum_{\mu < \nu}  e^{\mu} \wedge e^{\nu} \otimes (\langle s_{\mu}, s_{\nu} \rangle_h - \langle s_{\nu}, s_{\mu} \rangle_h) \\
            &= s^h \wedge s,
  \end{align*}
  where $s^h \mathpunct{:}= \sum_{\mu} e^{\mu} \otimes \langle s_{\mu}, \cdot \rangle_h  \in T^*M \otimes E^*$, $E^*$ is the complex dual bundle of $E$.

  In this special setting, the generalized Seiberg--Witten functional takes the form
  \begin{align}\label{hw}
    \mathcal{H}(A,\sigma) \mathpunct{:}= \int_M d\mathrm{vol} \left(|F_A - \sigma^h \wedge \sigma/2|^2 + |d_A \sigma|^2 + |d_A^* \sigma|^2 \right),
  \end{align}
  where $\sigma$ is a section of $T^*M \otimes E$.
  Equivalently, \eqref{hw} can be written as 
  \begin{equation}\label{eq:HW}
    \mathcal{H}(A,\sigma) = \int_M d\mathrm{vol} \left(|F_A|^2 + |\nabla_A \sigma|^2 + \langle \sigma, \mathrm{Ric}(\sigma) \rangle + |\sigma^h \wedge \sigma|^2/4 \right) d\mathrm{vol}.
  \end{equation}
\end{exmp}
The Euler--Lagrange equations \eqref{eq:A} and \eqref{eq:sigma} for the generalized Seiberg--Witten functional is transformed into the following explicit form for $\mathcal{H}$:
\begin{equation}\label{eq:EL-HW}
  \begin{cases}
    d^* F_A - \sum_{\mu} e^{\mu} \otimes \mathrm{Im}(\langle \nabla_{A, \mu} \sigma, \sigma \rangle_h) = 0,\\
    \nabla_A^*\nabla_A \sigma +   \sum_{\mu,\nu}e^{\mu} \otimes \mathrm{Ric}_{\mu\nu} \sigma^{\nu} 
    - \sum_{\mu} \mathrm{Im}(\langle \sigma, \sigma_{\mu} \rangle_h) \otimes \sigma_{\mu} = 0.
  \end{cases}
\end{equation}
In fact, in this special setting, we can write 
\[
\langle \nabla_A \sigma, B \sigma \rangle 
= \sum_{\mu} \langle \nabla_{A,\mu} \sigma, B_{\mu} \sigma \rangle = \sum_{\mu} \mathrm{Re}( \langle \nabla_{A,\mu} \sigma, \sigma \rangle_h B_{\mu}) 
= -\sum_{\mu} \mathrm{Im}(\langle \nabla_{A,\mu} \sigma, \sigma \rangle_h) \langle e^{\mu} \otimes i, B \rangle.
\] 
It follows that
\[
  \mathfrak{S}_A(\sigma)=-\sum_{\mu} \mathrm{Im}(\langle \nabla_{A, \mu} \sigma, \sigma \rangle_h) e^{\mu} \otimes i. 
\]
By \cref{Ric}, we  have
\[
  \mathrm{Ric}(\sigma) = \sum_{\mu,\nu}e^{\mu} \otimes \mathrm{Ric}_{\mu\nu} \sigma^{\nu},
\]
where $\mathrm{Ric}$ is the Ricci curvature of $M$, and
\begin{align*}
  \mathfrak{T}(\sigma) &= \sum_{\mu<\nu}  \langle (e_{\mu}e_{\nu} \otimes i \sigma, \sigma \rangle (e_{\mu}e_{\nu} \otimes i)\sigma \\
                       &= \sum_{\mu<\nu}  (\langle e_{\nu} \otimes i \sigma_{\mu}-e_{\mu} \otimes i \sigma_{\nu}, \sigma \rangle) (e_{\nu} \otimes i \sigma_{\mu}-e_{\mu} \otimes i \sigma_{\nu}) \\
                       &= \frac{1}{2}\sum_{\mu,\nu}  (\langle  i \sigma_{\mu}, \sigma_{\nu} \rangle - \langle  i \sigma_{\nu}, \sigma _{\mu} \rangle) (e_{\nu} \otimes i \sigma_{\mu}- e_{\mu} \otimes i \sigma_{\nu}) \\
                       &= \sum_{\mu}  (\langle  i \sigma_{\mu}, \sigma \rangle- \langle i \sigma, \sigma_{\mu} \rangle) \otimes i \sigma_{\mu} \\ 
                       &= -2 \sum_{\mu} \mathrm{Im}(\langle \sigma, \sigma_{\mu} \rangle_h) \otimes \sigma_{\mu}.
\end{align*}
\begin{rmk}\label{rmk: H}
  The $\mathrm{SU}(2)$ Kapustin--Witten functional can be interpreted as the $\mathrm{SL}(2,\mathbb{C})$ Yang--Mills functional with an additional gauge fixing term using the identification  $\mathfrak{sl}(2,\mathbb{C}) \cong \mathfrak{su}(2) \oplus i\mathfrak{su}(2)$\cite{Gagliardo-Uhlenbeck;Kapustin--Witten;Fixed;2012}. Likewise, \cref{hwf} can be interpreted as the $\mathrm{SU}(2)$ Yang--Mills functional with an additional gauge fixing term using the identification $\mathfrak{su}(2) \cong i \mathbb{R} \oplus \mathbb{C}$. More precisely, the diagonal embedding
  \begin{align*}
    \mathrm{U}(1) &\hookrightarrow \mathrm{SU}(2) \\
    \exp(i\varphi) &\mapsto 
    \begin{pmatrix}
      \exp(i\varphi) &0\\
      0 & \exp(-i\varphi)
    \end{pmatrix}
  \end{align*} 
  and the canonical action of $\mathrm{U}(1)$ on $\mathbb{C}$ induce the following identification
  \begin{align*}
    i \mathbb{R} \oplus \mathbb{C} &\rightarrow \mathfrak{su}(2) \\
    (ai,b+ic) &\mapsto a e_3 + b e_2 + c e_1,
  \end{align*}
  where
  \begin{align*}
    e_1 =\begin{pmatrix}
      0 & i\\
      i & 0
    \end{pmatrix},
    \quad
    e_2 =\begin{pmatrix}
      0 & 1\\
      -1 & 0
    \end{pmatrix},
    \quad
    e_3 =\begin{pmatrix}
      i & 0\\
      0 & -i
    \end{pmatrix}.
  \end{align*}
  The Lie bracket on $i \mathbb{R} \oplus \mathbb{C}$ is defined as
  \begin{align*}
    [(ai,\alpha),(bi,\beta)]=(-\mathrm{Im}(\overline{\alpha}\beta) i, ai \beta - bi \alpha).
  \end{align*}
  Moreover, one can check that the associative inner product on $i \mathbb{R} \oplus \mathbb{C}$ defined by $\langle (ai, \alpha),(bi, \beta) \rangle = ab + \mathrm{Re}(\overline{\alpha}\beta)$ agrees with the associative inner product on $\mathfrak{su}(2)$ defined by $\langle A, B \rangle = -\frac{1}{2}\mathrm{Tr}(AB)$ under this identification.

  Let $H$ be a principal $SU(2)$-bundle over $M$ with a reduction to the principal $U(1)$-bundle $P$ in \cref{hwf}. A connection $1$-form $\widetilde{A}$ over $H$ can be then decomposed as 
  \begin{align*}
    \widetilde{A} = A + \sigma,
  \end{align*}
  where $A$ is a connection $1$-form of $P$ and $\sigma$ is a section of the hermitian line bundle $E$ associated with $P$. Note that
  \begin{align*}
    F_{\widetilde{A}} = F_A +  d_A \sigma + \frac{1}{2}[\sigma,\sigma],
  \end{align*}
  and that
  \begin{align*}
    [\sigma,\sigma] = 2 \sum_{\mu < \nu} e^{\mu} \wedge e^{\nu} \otimes [\sigma_{\mu}, \sigma_{\nu}] 
    = -2\sum_{\mu < \nu} e^{\mu} \wedge e^{\nu} \otimes \mathrm{Im}(\overline{\sigma_{\mu}}\sigma_{\nu}) i = -\sigma^h \wedge \sigma.
  \end{align*}
  It follows that
  $
  |F_{\widetilde{A}}|^2  = |F_A - \frac{1}{2} \sigma^h \wedge \sigma|^2 + |d_A \sigma|^2.
  $ 
  Thus,
  \begin{align*}
    \mathcal{H}(A,\sigma) = \mathcal{YM}(\widetilde{A}) + \int_M d\mathrm{vol} \lvert d_A^* \sigma \rvert^2 := \int_M d\mathrm{vol} \left(|F_{\widetilde{A}}|^2 + \lvert d_A^* \sigma \rvert^2\right).
  \end{align*}
\end{rmk}

\section{Regularity of weak solutions of the generalized Seiberg--Witten equations}\label{sec:regularity}
We will establish in this section that a weak solution of $\mathcal{GSW}$ is indeed smooth. To do this, we will treat the equation for $A$ as a perturbation of the Yang--Mills equation, which is locally an elliptic equation for the connection 1-form in the local Coulomb gauge. The regularity of the classical Yang--Mills equation on $4$-manifolds can be found in \cite{Riviere;variations-Yang--Mills-Lagrangian;;2020}. By improving the regularity of $A$, we can bootstrap the regularity of $\sigma$, and by iterating this process, we will establish the smoothness of the weak solution. These regularity results will be crucial in proving the compactness and existence of weak solutions for $\mathcal{GSW}$ in \cref{sec:compactness}.
\begin{proof}[Proof of \cref{thm:regularity}]
  To do various estimates of the weak solution, let us write the equations of
  \eqref{eq:sigma} and \eqref{eq:A}  in terms of  some fixed smooth connection $\A$. Writing $A=\A+a$, then
  \begin{gather}
    \nabla_{\A} ^*\nabla _{\A}\sigma+\nabla _\A a\# \sigma+a\#\nabla _\A\sigma+a\#a\#\sigma+R_M\#\sigma+\sigma\#\sigma\#\sigma=0,\label{eq:sigma-Aref}\\
    d^*d a+\nabla _\A a\#a+a\#a\#a+\nabla _\A\sigma\#\sigma+a\#\sigma\#\sigma=0,\label{eq:A-Aref}
  \end{gather}
  where $\#$ denotes any multilinear map with smooth bounded coefficients.

  Now, let $A$ be in a Coulomb gauge on a neighbourhood $B$ of the point of consideration, that is $d^*a=0$ and 
  \begin{equation}
    \lVert a \rVert_{W^{1,2}(B)}\leq C\lVert F_A \rVert_{L^2(B)},
  \end{equation}
  we refer to \cite{Uhlenbeck;connections-curvature;;1982} for the local existence of Coulomb gauge.

  Firstly, we will bootstrap $a$. By the Bochner--Weitzenb\"ock formula, 
  \[
    d^*da+dd^*a=\Delta_{\mathrm{Hodge}} a=\nabla _\A^*\nabla _\A a+a\#\mathrm{Ric}+\nabla _\A a\# a+a\#a\#a,
  \]
  and noting that $A$ is in Coulomb gauge over $B$, \eqref{eq:A-Aref} yields 
  \begin{equation}\label{eq:a}
    \nabla_\A ^*\nabla_\A  a=a\#\nabla_\A a+a\#a\#a+\nabla _\A\sigma\#\sigma+a\#\sigma\#\sigma.
  \end{equation}
  The regularity of $A$ can be enhanced, employing the following growth estimates.
  \begin{clm}
    There exists $\alpha>0$, such that
    \begin{equation}\label{eq:growth-a}
      \sup_{\substack{x_0\in B_{1/2}\\0<\rho<1/4}}\rho^{-4\alpha}\int_{B_\rho(x_0)}\lvert a \rvert^4dx<+\infty.
    \end{equation}
  \end{clm}
  We follow the argument of
  \cite{Riviere;variations-Yang--Mills-Lagrangian;;2020}. In fact, let
  $\epsilon>0$ to be fixed later, at this moment we assume $\epsilon<1$. There exists $\rho_0>0$ such that
  \begin{equation}\label{eq:assumption-a}
    \sup_{\substack{x_0\in B_{1/2}\\0<\rho<\rho_0}}\int_{B_\rho(x_0)}\left( \lvert a \rvert^4+\lvert \sigma \rvert^4+\lvert \nabla_\A a \rvert^2+\lvert \nabla _\A\sigma \rvert^2 \right)dx<\epsilon,
  \end{equation}
  because $\lVert a \rVert_{W^{1,2}(B)}$, $\lVert \nabla_{\A}\sigma \rVert_{L^2(B)}$ and $\lVert \sigma \rVert_{L^4(B)}$ are all bounded.
  For any $x_0\in B_{1/2}$ and any $0<\rho_1<\rho_0$, consider the following equation
  \[
    \begin{cases}
      \Delta_{\mathrm{rough}} \tilde a\mathpunct{:}=\nabla_\A ^*\nabla_\A a=a\#\nabla_\A a+a\#a\#a+\nabla_\A\sigma\#\sigma+a\#\sigma\#\sigma,&x\in B_{\rho_1}(x_0)\\
      \tilde a=0,&x\in\partial B_{\rho_1}(x_0).
    \end{cases}
  \]
  The classical elliptic estimates show, for some constant $C$ independent of $\rho_1$, 
  \begin{align*}
    \lVert \tilde a \rVert_{L^4(B_{\rho_1}(x_0))}
    &\leq C\lVert \tilde a \rVert_{W^{2,4/3}(B_{\rho_1}(x_0))}\\
    &\leq C\lVert a\#\nabla_\A a+a\#a\#a+\nabla _\A\sigma\#\sigma+a\#\sigma\#\sigma \rVert_{L^{4/3}(B_{\rho_1}(x_0))}\\
    &\leq C\left( \lVert a \rVert_{L^4(B_{\rho_1}(x_0))}\lVert \nabla_\A a \rVert_{L^2(B_{\rho_1}(x_0))}+\lVert \nabla _\A\sigma \rVert_{L^2(B_{\rho_1}(x_0))}\lVert \sigma \rVert_{L^4(B_{\rho_1}(x_0))} \right)\\
    &\qquad +C\lVert a \rVert_{L^4(B_{\rho_1}(x_0))}\left(\lVert a \rVert_{L^4(B_{\rho_1}(x_0))}^2+\lVert \sigma \rVert_{L^4(B_{\rho_1}(x_0))}^2\right)\\
    &\leq C\left( \lVert \nabla_\A a \rVert_{L^2(B_{\rho_1}(x_0))}+\lVert a \rVert_{L^4(B_{\rho_1}(x_0))}^2+\lVert \sigma \rVert_{L^4(B_{\rho_1}(x_0))}^2 \right)\lVert a \rVert_{L^4(B_{\rho_1}(x_0))}\\
    &\qquad+C \lVert \nabla _\A\sigma \rVert_{L^2(B_{\rho_1}(x_0))}\lVert \sigma \rVert_{L^4(B_{\rho_1}(x_0))}.
  \end{align*}
  Furthermore, since $\lVert \sigma \rVert_{L^\infty(B_{\rho_1}(x_0))}$ is bounded by \cref{cor:L-infty-boundedness-via-Morse}, we can use the assumption \eqref{eq:assumption-a} to obtain
  \begin{equation}\label{eq:L4-a}
    \begin{split}
      &\lVert \tilde a \rVert_{L^4(B_{\rho_1}(x_0))}\\
      &\qquad\leq C\left( \lVert \nabla_\A a \rVert_{L^2(B_{\rho_1}(x_0))}+\lVert a \rVert_{L^4(B_{\rho_1}(x_0))}^2+\lVert \sigma \rVert_{L^4(B_{\rho_1}(x_0))}^2 \right)\lVert a \rVert_{L^4(B_{\rho_1}(x_0))}\\
      &\qquad\qquad+C\lVert \nabla _\A\sigma \rVert_{L^2(B_{\rho_1}(x_0))}\rho_1\\
      &\qquad\leq C\epsilon\left( \lVert a \rVert_{L^4(B_{\rho_1}(x_0)}+\rho_1 \right).
    \end{split}
  \end{equation}

  Now, note that the difference $\bar a=a-\tilde a$ is harmonic on $B_{\rho_1}(x_0)$, hence $\lvert \bar{a} \rvert^4$ is sub-harmonic, and 
  \[
    \frac{d}{dr}\left[ \frac{1}{r^4}\int_{B_r(x_0)}\lvert \bar a \rvert^4(x)dx \right]\geq0.
  \]
  In particular,
  \[
    \int_{B_{\frac{\rho_1}{4}}(x_0)}\lvert \bar a \rvert^4\leq 2^{-8}\int_{B_{\rho_1}(x_0)}\lvert \bar a \rvert^4,
  \]
  and
  \begin{align*}
    \int_{B_{\frac{\rho_1}{4}}(x_0)}\lvert a \rvert^4&\leq 8\int_{B_{\frac{\rho_1}{4}}(x_0)}\left( \lvert \bar a \rvert^4+\lvert \tilde a \rvert^4 \right)
    \leq 2^{-5}\int_{B_{\rho_1}(x_0)}\lvert \bar a \rvert^4+8\int_{B_{\frac{\rho_1}{2}}(x_0)}\lvert \tilde a \rvert^4\\
                                                     &\leq 2^{-2}\int_{B_{\rho_1}(x_0)}\lvert a \rvert^4+16\int_{B_{\frac{\rho_1}{2}}(x_0)}\lvert \tilde a \rvert^4.
  \end{align*}
  By the $L^4$-estimates of $\tilde a$, see \eqref{eq:L4-a}, we have
  \[
    \int_{B_{\frac{\rho_1}{4}}(x_0)}\lvert a \rvert^4\leq \left( 2^{-2}+C\epsilon^4 \right)\lVert a \rVert_{L^4(B_{\rho_1}(x_0))}^4+C\epsilon^4 \rho_1^4.
  \]
  If we choose $\epsilon$ such that $C\epsilon^4<2^{-2}$, then we conclude that
  \[
    \int_{B_{\frac{\rho_1}{4}}(x_0)}\lvert a \rvert^4\leq \frac{1}{2}\int_{B_{\rho_1}(x_0)}\lvert a \rvert^4+2^{-2}\rho_1^4.
  \]
  Iteration implies
  \[
    \int_{B_{4^{-k}\rho_1}(x_0)}\lvert a \rvert^4\leq 2^{-k}\int_{B_{\rho_1}(x_0)}\lvert a \rvert^4+\sum_{i=1}^k 2^{-(i-1)}(2^{-2}(2^{-2(k-i)}\rho_1)^4)
    \leq 2^{-k}\int_{B_{\rho_1}(x_0)}\lvert a \rvert^4+2^{-8k}\rho_1^4,
  \]
  thus
  \[
    \int_{B_{4^{-k}\rho_1}(x_0)}\lvert a \rvert^4\leq 2^{-k}\left( \int_{B_{\rho_1}(x_0)}\lvert a \rvert^4+\rho_1^4 \right),
  \]
  i.e.,
  \[
    \frac{1}{(2^{-2k}\rho_1)^{4\alpha}}\int_{B_{2^{-2k}\rho_1}(x_0)}\lvert a \rvert^4\leq 2^{k(8\alpha-1)}\rho_1^{-4\alpha}\left( \rho_1^{4\alpha}\cdot \frac{1}{\rho_1^{4\alpha}}\int_{B_{\rho_1}(x_0)}\lvert a \rvert^4+\rho_1^4  \right),
  \]
  thus if we take $\alpha<1/8$, then 
  \[
    \frac{1}{(2^{-2k}\rho_1)^{4\alpha}}\int_{B_{2^{-2k}\rho_1}(x_0)}\lvert a \rvert^4
    \leq \frac{1}{\rho_1^{4\alpha}}\int_{B_{\rho_1}(x_0)}\lvert a \rvert^4+\rho_1^{4-4\alpha}.
  \]
  Now, assume $2^{-2k_0}<\rho_0\leq 2^{-2(k_0-1)}$, and for any $0<\rho<1/4$, 
  \begin{itemize}
    \item if $2^{-2(k+k_0+1)}\leq\rho<2^{-2(k+k_0)}$, then $2^{-4}\rho_0\leq\rho_1=2^{2k}\rho<2^{-2k_0}<\rho_0$, and
      \begin{align*}
        \frac{1}{\rho^{4\alpha}}\int_{B_\rho(x_0)}\lvert a \rvert^4
  &=\frac{1}{(2^{-2k}\rho_1)^{4\alpha}}\int_{B_{2^{-2k}\rho_1}(x_0)}\lvert a \rvert^4
  \leq \frac{1}{\rho_1^{4\alpha}}\int_{B_{\rho_1}(x_0)}\lvert a \rvert^4+\rho_1^{4-4\alpha}\\
  &\leq \frac{1}{(2^{-4}\rho_0)^{4\alpha}}\int_{B_{\rho_0}(x_0)}\lvert a \rvert^4+\frac{1}{(2^{-4}\rho_0)^{4\alpha-4}}
  \leq C(\rho_0)\left( \int_{B_{\rho_0}(x)}\lvert a \rvert^4+1 \right);
      \end{align*}
    \item if $\rho\geq2^{-2k_0}\geq 2^{-2}\rho_0$, then
      \[
        \frac{1}{\rho^{4\alpha}}\int_{B_\rho(x_0)}\lvert a \rvert^4\leq \frac{1}{(2^{-2}\rho_0)^{4\alpha}}\int_{B_{1/4}(x_0)}\lvert a \rvert^4,
      \]
  \end{itemize}
  In either case, $\rho^{-4\alpha}\int_{B_\rho(x_0)}\lvert a \rvert^4$ is bounded. Hence, the Claim is proved.

  Similar to \eqref{eq:L4-a}, the elliptic estimates of \eqref{eq:a} imply
  \[
    \lVert \nabla_\A ^*\nabla_\A a \rVert_{L^{4/3}(B_\rho(x_0))}\leq C(\lVert a \rVert_{L^4(B_\rho(x_0))}+\rho),
  \]
  where $C$ is a constant, depending on the energy $\mathcal{GSW}(A,\sigma)$ and $\lVert \sigma \rVert_{L^\infty}$. Thus, 
  \begin{align*}
    \sup_{\substack{x_0\in B_{1/2}\\ 0<\rho<1/4}} \frac{1}{\rho^{4\alpha/3}}\int_{B_\rho(x_0)}\lvert \nabla_\A ^*\nabla_\A a \rvert^{4/3}
  &\leq C \sup_{\substack{x_0\in B_{1/2}\\ 0<\rho<1/4}} \frac{1}{\rho^{4\alpha/3}} \left[ \left( \int_{B_\rho(x_0)}\lvert a \rvert^4 \right)^{1/3}+\rho^{4/3} \right]\\
  &\leq C \sup_{\substack{x_0\in B_{1/2}\\ 0<\rho<1/4}}  \left[ \left( \frac{1}{\rho^{4\alpha}}\int_{B_\rho(x_0)}\lvert a \rvert^4 \right)^{1/3}+\rho^{4/3(1-\alpha)} \right].
  \end{align*}
  which is bounded by the Claim. That is, $a\in L^{3-\alpha,4/3}(B_{1/2})$, where the Morrey space $L^{\lambda,p}(B)$ of $\mathfrak{g}_E$-valued 1-forms on $B$ is the completion of $\Omega^1(B,\mathfrak{g}_E)$ with respect to the following norm
  \[
    \lVert a \rVert_{L^{\lambda,p}(B)}\mathpunct{:}=\lVert a \rVert_{L^p(B)}+\sup_{\substack{x\in B\\0<\rho<1/4}}\left( \frac{1}{\rho^{n-\lambda p}}\int_{B_\rho(x)\cap B}\lvert a \rvert^p \right)^{1/p}.
  \]
  By standard elliptic regularity and by Adams-Sobolev embeddings \cite{Adams;Riesz-potentials;Duke;1975}, we conclude that $\nabla_\A a\in L^p_{\mathrm{loc}}(B_{1/2})$ for some $p>2$. 

  Finally, we can bootstrap to obtain the desired regularity. In fact, since
  $\nabla_\A a\in L^p_{\mathrm{loc}}(B_{1/2})$, we know that the RHS of
  \eqref{eq:a} belongs to some $L^{8/5}_{\mathrm{loc}}(B_{1/2})$. The standard
  elliptic estimates imply that $a\in W^{2,8/5}_{\mathrm{loc}}(B_{1/2})$ and the
  Sobolev embedding implies that $\nabla_\A a\in
  L^{8/3}_{\mathrm{loc}}(B_{1/2})$. This implies the RHS of \eqref{eq:a} belongs
  to $L^2_{\mathrm{loc}}(B_{1/2})$, and $a\in W^{2,2}_{\mathrm{loc}}(B_{1/2})$
  by standard elliptic estimates. Now, the RHS of \eqref{eq:sigma-Aref} belongs
  to $L^q_{\mathrm{loc}}(B_{1/2})$ for any $1<q<2$, thus the elliptic estimates
  implies $\sigma\in W^{2,q}_{\mathrm{loc}}(B_{1/2})$. The Sobolev embedding
  implies that $\nabla_\A\sigma\in L^{q'}_{\mathrm{loc}}(B_{1/2})$ for any
  $1<q'<4$. Then the RHS of \eqref{eq:a} belongs to
  $L^{q''}_{\mathrm{loc}}(B_{1/2})$ for any $1<q''<4$, thus $a\in
  W^{2,q''}_{\mathrm{loc}}(B_{1/2})$ and is H\"older continuous. This also
  improves the RHS of \eqref{eq:sigma-Aref} to $L^{q'}_{\mathrm{loc}}(B_{1/2})$,
  and $\sigma\in W^{2,q'}_{\mathrm{loc}}(B_{1/2})$, where $1<q'<4$, thus
  $\sigma$ is also H\"older continuous. The routine bootstrap argument then
  implies the smoothness of $(A,\sigma)$.
\end{proof}
\begin{rmk}
In fact, a glueing argument similar to \cite{Jost-Kessler-Wu-Zhu;Geometric}*{Theorem~3.2} will show that the local gauge can be pieced together to obtain a global $W^{2,2}$ gauge, such that $g(A,\sigma)$ is smooth.
\end{rmk}

\section{Compactness and existence theorems for the \texorpdfstring{$\mathcal{H}$}{H}-functional}\label{sec:compactness}
In this section, we will complete the proof of the existence of \cref{thm:main} by verifying the Palais--Smale compactness condition. As a preparation, we will first show the coerciveness of $\mathcal{H}$, where the abelian structure group and the closed sub-manifold $\mathcal{V}$ play an important role as explained in \cref{sec:intro}. This coerciveness provides a natural boundedness of Palais--Smale sequences in $\mathcal{V}$. Then, the verification of the Palais--Smale condition is somewhat routine; however, we still need to restrict to the closed sub-manifold $\mathcal{V}$ to prove the strong convergence of connections $A_n$, as seen in \cref{step:clm,step:An-convergence}. It should also be noted that the strong convergence of $\sigma_n$ depends on the strong convergence of $A_n$, as seen in the process of term $\mathrm{II}$ in \cref{step:sigman-convergence}.

Firstly, we will choose a suitable Lie group as the abelian gauge group, similar to the Seiberg--Witten case, which will act smoothly on $\mathscr{A}^{1,2}\times W^{1,2}(\Omega^1(E))$. Let $\mathscr{G}_0=\exp\left( \sqrt{-1}W^{2,2}(M,\mathbb{R}) \right)$, which is a Lie group and can be identified through the exponential map with the quotient $W^{2,2}(M,\mathbb{R})$ under the relation $g_1\sim g_2$ if and only if $g_1-g_2\in 2\pi \mathbb{Z}$ for almost all $x\in M$. The gauge group is then defined as $\mathscr{G}^{2,2}=\sqcup g\cdot \mathscr{G}_0$, where the union is taken over all components of $C^\infty(M,S^1)$. $\mathscr{G}^{2,2}$ is a Lie group and acts smoothly on $\mathscr{A}^{1,2}\times W^{1,2}(\Omega^1(E))$, as shown in \cite{Jost-Peng-Wang;variational;calc;1996}.

In addition to introducing the subspace $\mathcal{V}$, we will also introduce the following functions for $(A,\sigma)\in \mathcal{V}$. It is important to note that $\tau$ exhibits a quadratic dependence on $\sigma$. Let
\[
  s(x)\mathpunct{:}=\min_{\lvert \sigma(x) \rvert=1}\frac{\left\langle \sigma(x),\mathrm{Ric}(\sigma)(x) \right\rangle}{\lvert \tau(\sigma)(x) \rvert},\quad
  s_0 \mathpunct{:}=\min_{x\in M}s(x).
\]
\begin{rmk}
    In the classical Seiberg--Witten case, the function $s(x)$ is exactly the scalar curvature.
\end{rmk}
\begin{lem}\label{lem:coercive}
  The functional $\mathcal{H}$ is coercive on $\mathcal{V}$, i.e., there exists a constant $c>0$ such that for each $(A,\sigma)\in \mathcal{V}$, we have 
  \[
    \mathcal{H}(A,\sigma)\geq c^{-1}\left( \lVert g^{-1}\sigma \rVert_{W^{1,2}}+\lVert g(A) \rVert_{W^{1,2}}-1 \right),
  \]
  for some gauge transformation $g\in \mathscr{G}^{2,2}$. 
\end{lem}
\begin{proof}
  Firstly, we estimate the term containing $g(A)$, which can be done as in \cite{Jost-Peng-Wang;variational;calc;1996}*{Lem.~2.5}. In fact, we can first fix the gauge by considering $\zeta\in W^{2,2}(M , \mathbb{R})$ and solving
  \[
    d^*d\zeta=\sqrt{-1} d^*(A-\A)\in L^2(\Omega^1(\ad P_{G/\mathbb{Z}_2})),
  \]
  then $g_0=\exp(\sqrt{-1}\zeta)\in \mathscr{G}_0$ satisfies 
  \[
    d^*\left( g_0(A)-\A \right)=d^*\left( g_0^{-1}dg_0+A-\A \right)=0.
  \]

  Since the component group of $C^\infty(M,S^1)$ is isomorphic to $H^1(M,\mathbb{Z})$, by a well-known result about harmonic maps from $M$ to $S^1$, see \cite{Eells-Lemaire;Another-report;;1988}, for any component of $C^\infty(M,S^1)$, there exists a unique map $g_1\in C^\infty(M,S^1)$ satisfying 
  \[
    d^*(g_1^{-1}dg_1)=0,\quad g_1(x_0)=1,
  \]
  where $x_0\in M$ is a fixed point. In fact, the harmonic maps into $S^1$ solve
  \[
    \Delta g_1+g_1\lvert \nabla  g_1 \rvert^2=0\implies d^*\left( g_1^{-1}dg_1 \right)=\lvert dg_1 \rvert^2+g_1^{-1}\Delta g_1=g_1^{-1}\left( g_1\lvert \nabla  g_1 \rvert^2+\Delta g_1 \right)=0,
  \]
  where $\Delta$ is the geometric Laplace--Beltrami operator, which differs from the analytic Laplacian by a minus sign.
  By the Hodge decomposition, we have for $g = g_1g_0$,
  \[
    \lVert g(A)-\A \rVert_{W^{1,2}}\leq c\left( \lVert d^*\left( g(A)-\A \right) \rVert_{L^2}+\lVert d\left( g(A)-\A \right) \rVert_{L^2}+\lVert H\left( g(A)-\A \right) \rVert_{L^2} \right),
  \]
  where $H(\cdot)$ is the harmonic part. Note that, 
  \[
    g(A)-\A=g_1^{-1}dg_1+g_0(A)-\A\implies d^*\left( g(A)-\A \right)=0.
  \]

  Furthermore, considering that the harmonic component of $g_1^{-1}dg_1$ lies in $H^1(M,\mathbb{Z})$, and the Jacobi torus $H^1(M,\mathbb{R})/H^1(M,\mathbb{Z})$ is compact, it is possible to select a component of $C^\infty(M,S^1)$ in such a way that the boundedness of the harmonic part of $g_1^{-1}dg_1$ is ensured. The harmonic part of $g(A)-\A$ encompasses both the harmonic portion of $g_0(A)-\A$ and the harmonic part of $g_1^{-1}dg_1$, which is $g_1^{-1}dg_1$ itself, and also bounded. In conclusion,
  \begin{equation}\label{eq:gA-W12}
    \lVert g(A)-\A \rVert_{W^{1,2}}\leq c\lVert d(g(A)-\A) \rVert_{L^2}+c
    \leq c\lVert d(g(A)) \rVert_{L^2}+c\leq c\lVert F_A \rVert_{L^2}+c,
  \end{equation}
  where we use the fact that $G$ is abelian in the last inequality.

  Next, we estimate $\lVert g^{-1}\sigma \rVert_{W^{1,2}}$. Note first that, by the assumption $\lvert \sigma \rvert^2\leq \lambda_0\lvert \tau(\sigma) \rvert$, we have 
  \[
    \lVert g^{-1}\sigma \rVert_{L^2}=\lVert \sigma \rVert_{L^2}\leq c \lVert  \sigma  \rVert_{L^4}\leq c\lVert \tau(\sigma) \rVert_{L^2}^{1/2}\leq c\lVert \tau(\sigma) \rVert_{L^2}+c,
  \]
  and by \eqref{eq:gA-W12},
  \begin{align*}
    \lVert \nabla_\A (g^{-1}\sigma) \rVert_{L^2}&\leq\lVert \nabla _{g(A)}(g^{-1}\sigma)-(g(A)-\A)(g^{-1}\sigma) \rVert_{L^2}\\
                                                &\leq c \left( \lVert \nabla _{g(A)}(g^{-1}\sigma) \rVert_{L^2}+\lVert g(A)-\A \rVert_{L^4}\lVert g^{-1}\sigma \rVert_{L^4} \right)\\
                                                &\leq c\left( \lVert g^{-1}\nabla _{A} (gg^{-1}\sigma) \rVert_{L^2}+\lVert F_A \rVert_{L^2}^2+\frac{1}{8}\lVert \tau(\sigma) \rVert_{L^2} \right)+c\\
                                                &\leq c\left( \lVert \nabla _{A} \sigma \rVert_{L^2}+\lVert F_A \rVert_{L^2}^2+\frac{1}{8}\lVert \tau(\sigma) \rVert_{L^2} \right)+c.
  \end{align*}

  Finally, since $M$ is compact, we may assume that $\lVert s \rVert_{L^\infty(M)}^2\leq c$, and
  \begin{align*}
    \mathcal{H}(A,\sigma)&=\int_M \left( \lvert F_A \rvert^2+\lvert \nabla _A\sigma \rvert^2+\left\langle \sigma,\mathrm{Ric}(\sigma) \right\rangle+\lvert \tau(\sigma) \rvert^2/4 \right)\\
                         &\geq \int_M \left( \lvert F_A \rvert^2+\lvert \nabla _A\sigma \rvert^2+\frac{1}{8}\lvert \tau(\sigma) \rvert^2-2s^2 \right)\\
                         &\qquad+\frac{1}{8}\int_M\left( 16s^2+8s\lvert \tau(\sigma) \rvert+\lvert \tau(\sigma) \rvert^2 \right)\\
                         &\geq\int_M\left( \lvert F_A \rvert^2+\lvert \nabla _A\sigma \rvert^2+\frac{1}{8}\lvert \tau(\sigma) \rvert^2 \right) -c.
  \end{align*}

  In summary, we have proved there exists a $g\in \mathscr{G}^{2,2}$, such that 
  \[
    \lVert g^{-1}\sigma \rVert_{W^{1,2}}+\lVert g(A) \rVert_{W^{1,2}}\leq c \mathcal{H}(A,\sigma)+c.
  \]
\end{proof}
Now, we are ready to verify the Palais--Smale compactness condition.
\begin{proof}[Proof of \cref{thm:main}]
  By Lemma~\autoref{lem:coercive} we know that there exist $g_n\in \mathscr{G}^{2,2}$, such that 
  \[
    \lVert g_n(A_n) \rVert_{W^{1,2}}+\lVert g_n^{-1}\sigma_n \rVert_{W^{1,2}}\leq c,
  \]
  where $c$ is a uniform constant independent of $n$. To simplify our notation, we will use $A_n$ and $\sigma_n$ to represent $g_n(A_n)$ and $g_n^{-1}\sigma_n$, respectively. By the compact embedding theorem, we can conclude that, for a subsequence, we have 
  \begin{enumerate}
    \item $A_n \rightharpoonup A_\infty$ weakly in $W^{1,2}(\Omega^1(\ad P_{G/\mathbb{Z}_2}))$, and $\sigma_n \rightharpoonup \sigma_\infty$ weakly in $W^{1,2}(\Omega^1(E))$;
    \item $A_n \rightharpoonup A_\infty$ weakly in $L^{4}(\Omega^1(\ad P_{G/\mathbb{Z}_2}))$, and $\sigma_n \rightharpoonup \sigma_\infty$ weakly in $L^{4}(\Omega^1(E))$;
    \item for any $1<p<4$, $A_n \to A_\infty$ strongly in $L^{p}(\Omega^1(\ad P_{G/\mathbb{Z}_2}))$, and $\sigma_n \to \sigma_\infty$ strongly in $L^{p}(\Omega^1(E))$.
  \end{enumerate}

  \step We first verify that $(A_\infty,\sigma_\infty)$ is a weak solution of \eqref{eq:EL-HW}. In fact, for any 1-form $a\in W^{1,2}(\Omega^1(\ad P_{G/\mathbb{Z}_2}))$, by the first Palais--Smale condition, we have,
  \[
    \int_M\left\langle F_{A_n},da \right\rangle+\left\langle \nabla _{A_n}\sigma_n, a\sigma_n \right\rangle=d \mathcal{H}(A_n,\sigma_n)(a)=o(1).
  \]
  Note that 
  \[
    \int_{M}\left\langle F_{A_n}, da \right\rangle-\int_M\left\langle F_{A_\infty}, da \right\rangle= \int_M \left\langle F_{A_n}-F_{A_\infty}, da \right\rangle=o(1),
  \]
  because $F_{A_n}\rightharpoonup F_{A_\infty}$ in $L^2(\Omega^1(\ad P_{G/\mathbb{Z}_2}))$ and $da\in L^2(\Omega^1(\ad P_{G/\mathbb{Z}_2}))$. Similarly, if we write $a_n=A-\A$, and $a_\infty=A_\infty-\A$, then
  \begin{align*}
    &\int_M\left\langle \nabla _{A_n}\sigma_n,a\sigma_n \right\rangle-\int_M\left\langle \nabla _{A_\infty}\sigma_\infty,a\sigma_\infty \right\rangle\\
    &\qquad=\int_M\left\langle \nabla _\A\sigma_n+a_n\sigma_n,a\sigma_n \right\rangle-\int_M\left\langle \nabla _\A\sigma_\infty+a_\infty\sigma_\infty,a\sigma_\infty \right\rangle\\
    &\qquad=\int_M\left\langle \nabla _\A(\sigma_n-\sigma_\infty)+a_n(\sigma_n-\sigma_\infty), a\sigma_n \right\rangle-\int_M\left\langle \nabla _\A\sigma_\infty+a_n\sigma_\infty,a(\sigma_\infty-\sigma_n) \right\rangle,
  \end{align*}
  which tends to zero by the weak convergence of $(A_n,\sigma_n)$. Thus, we have proved that 
  \[
    \int_M\left\langle F_{A_\infty},da \right\rangle+\left\langle \nabla _{A_\infty}\sigma_\infty,a\sigma_\infty \right\rangle=0,\quad \forall a\in W^{1,2}(\Omega^1(\ad P_{G/\mathbb{Z}_2})),
  \]
  i.e., $(A_\infty,\sigma_\infty)$ solves the first equation of \eqref{eq:EL-HW} weakly. An analogous argument shows that $(A_\infty,\sigma_\infty)$ also solves the second equation weakly. Therefore, $(A_\infty,\sigma_\infty)$ is a weak solution of \eqref{eq:EL-HW}. By the regularity \cref{thm:regularity}, we know that $(A_\infty,\phi_\infty)$ is smooth.

  \step \label{step:clm} In order to prove the strong convergence of $A_n\to A_\infty$ in $W^{1,2}(\Omega^1(\ad P_{G/\mathbb{Z}_2}))$, we need to make the following claim. The technical contribution is due to Taubes \cite{Taubes;equivalence;Comm;1980}, and the same technique is used to prove the $L^\infty$-boundedness of $\sigma$ in \cite{Jost-Peng-Wang;variational;calc;1996}*{Lemma~3.2} for the Seiberg--Witten functional. We can also prove the $L^\infty$-boundedness of $\sigma$ for the generalized Seiberg--Witten functional based on this technique, provided that the weak solution is in $\mathcal{V}$. However, our argument in \cref{thm:regularity} does not have any restrictive conditions, unlike the one in the definition of $\mathcal{V}$.

  \begin{clm}
    Let $c_0 \mathpunct{:}=\sqrt{-2s_0\lambda_0}$ when $s_0<0$ and $c_0 \mathpunct{:}=1$ when $s_0\geq0$. Then, for 
    \[
      \Omega_n \mathpunct{:}=\left\{ x\in M: f(\sigma_n) \mathpunct{:}=\lvert \sigma_n \rvert-c_0>0 \right\}, \quad
      \nu_n \mathpunct{:}=\sigma_n/\lvert \sigma_n \rvert, \quad \sigma_n\in \Omega_n,
    \]
    we have, as $n\to\infty$  
    \[
      \int_{\Omega_n}\lvert \left\langle \nabla _{A_n}\sigma_n,\nu_n \right\rangle \rvert^2\to0.
    \]
  \end{clm}
  In fact, note first that, for
  $$
  \eta_n=
  \begin{cases}
    f(\sigma_n)\nu_n,&\lvert \sigma_n \rvert>c_0,\\
    0,&\lvert \sigma_n \rvert\leq c_0,
  \end{cases}
  $$
  by the assumption that $(A_n,\sigma_n)\in \mathcal{V}$, we have
  \begin{align*}
    d \mathcal{H}(A_n,\sigma_n)(\eta_n)&=\int_M\left\langle \nabla _{A_n}\sigma_n,\nabla _{A_n}\eta_n \right\rangle+\left\langle \mathrm{Ric}(\sigma_n),\eta_n \right\rangle+\frac{1}{2}\left\langle \mathfrak{T}(\sigma_n),\eta_n \right\rangle\\
                                       &=\int_{\Omega_n}\left\langle \nabla _{A_n}\sigma_n,\nabla _{A_n}\eta_n \right\rangle+\frac{f(\sigma_n)}{\lvert \sigma_n \rvert}\left( \left\langle \mathrm{Ric}(\sigma_n),\sigma_n \right\rangle+\lvert \tau(\sigma_n) \rvert^2/2 \right)\\
                                       &\geq\int_{\Omega_n}\left\langle \nabla _{A_n}\sigma_n,\nabla _{A_n}\eta_n \right\rangle+\frac{f(\sigma_n)\lvert \tau(\sigma_n) \rvert}{2\lvert \sigma_n \rvert}\left( 2s_0+\lvert \tau(\sigma_n) \rvert  \right)\\
                                       &\geq\int_{\Omega_n}\left\langle \nabla _{A_n}\sigma_n,\nabla _{A_n}\eta_n \right\rangle,
  \end{align*}
  in either case, since $f(\sigma_n)\geq0$ and $\lvert \sigma_n \rvert^2+2s_0\lambda_0\geq0$ on $\Omega_n$.

  We can continue the computation on $\Omega_n$ as follows: since $d\lvert \sigma_n \rvert=\left\langle \nabla _{A_n}\sigma_n,\nu_n \right\rangle$ on $\Omega_n$,
  \begin{align*}
    \nabla_{A_n}\sigma_n&=\nabla_{A_n}(\lvert \sigma_n \rvert\nu_n)=\lvert \sigma_n \rvert\nabla_{A_n}\nu_n+d\lvert \sigma_n \rvert\nu_n\\
                        &=\lvert \sigma_n \rvert\nabla_{A_n}\nu_n+\left\langle \nabla_{A_n}\sigma_n,\nu_n \right\rangle\nu_n\\
                        &=\left( f(\sigma_n)\nabla_{A_n}\nu_n+c_0\nabla_{A_n}\nu_n \right)+\left\langle \nabla_{A_n}\sigma_n,\nu_n \right\rangle \nu_n,
  \end{align*}
  we have,
  \begin{align*}
    \langle \nabla_{A_n} \sigma_n, \nabla_{A_n} \eta_n \rangle 
      &= \left\langle \nabla_{A_n}\sigma_n,d\lvert \sigma_n \rvert\nu_n \right\rangle+\left\langle \nabla_{A_n}\sigma_n,f(\sigma_n)\nabla _{A_n}\nu_n \right\rangle\\
      &=\left\langle \nabla_{A_n}\sigma_n,\nu_n \right\rangle^2+\left\langle f(\sigma_n)\nabla_{A_n}\nu_n, \nabla_{A_n}\sigma_n \right\rangle\\
      &=\left\langle \nabla_{A_n}\sigma_n,\nu_n \right\rangle^2+\left\langle f(\sigma_n)\nabla_{A_n}\nu_n, f(\sigma_n)\nabla_{A_n}\nu_n+c_0\nabla_{A_n}\nu_n+\left\langle \nabla_{A_n}\sigma_n,\nu_n \right\rangle\nu_n \right\rangle\\
      &\geq c_0f(\sigma_n)\lvert \nabla_{A_n}\nu_n \rvert^2+\left\langle \nabla_{A_n}\sigma_n,\nu_n \right\rangle^2+f^2(\sigma_n) \lvert \nabla_{A_n}\nu_n \rvert^2\\
      &\qquad-f(\sigma_n)\lvert \left\langle \nabla_{A_n}\sigma_n,\nu_n \right\rangle \rvert\lvert \nabla_{A_n}\nu_n \rvert\\
      &\geq c_0f(\sigma_n)\lvert \nabla_{A_n}\nu_n \rvert^2+\frac{1}{2}\left\langle \nabla_{A_n}\sigma_n,\nu_n \right\rangle^2.
  \end{align*}
  Therefore,
  \begin{multline*}
    d \mathcal{H}(A_n,\sigma_n)(\eta_n)
    \geq\int_{\Omega_n}\left\langle \nabla _{A_n}\sigma_n,\nabla _{A_n}\eta_n \right\rangle\\
    \geq\int_{\Omega_n}\frac{1}{2}\left\langle \nabla _{A_n}\sigma_n,\nu_n \right\rangle^2+c_0f(\sigma_n)\lvert \nabla _{A_n}\nu_n \rvert^2
    \geq \frac{1}{2}\int_{\Omega_n}\left\langle \nabla _{A_n}\sigma_n, \nu_n \right\rangle^2,
  \end{multline*}
  that is
  \[
    \int_{\Omega_n}\lvert\left\langle  \nabla _{A_n}\sigma_n,\nu_n  \right\rangle\rvert^2\leq 2 d \mathcal{H}(A_n,\sigma_n)(\eta_n)\leq 2\lVert d \mathcal{H}(A_n,\sigma_n) \rVert_{W^{-1,2}}\lVert \eta_n \rVert_{W^{1,2}}.
  \]
  The claim follows provided $\lVert \eta_n \rVert_{W^{1,2}}\leq c$. In fact, 
  \[
    \lVert \eta_n \rVert_{L^2}^2= \int_{\Omega_n} \left( \lvert \sigma_n \rvert-c_0 \right) ^2 \leq \int_{M}\lvert \sigma_n \rvert^2+c 
    \leq \lVert \sigma_n \rVert_{W^{1,2}}^2+c;
  \]
  and, since $\nabla_\A \eta_n=d\lvert \sigma_n \rvert\nu_n+f(\sigma_n)\nabla_\A \nu_n$ on $\Omega_n$,
  \begin{align*}
    \lVert \nabla_\A \eta_n \rVert_{L^2}^2&=\lVert d\lvert \sigma_n \rvert\nu_n+f(\sigma_n)\nabla_\A \nu_n \rVert_{L^2(\Omega_n)}^2\\
                                          &\leq c\left( \lVert \nabla_\A \sigma_n \rVert_{L^2}^2+\left\lVert \frac{\lvert \sigma_n \rvert-c_0}{\lvert \sigma_n \rvert} \right\rVert_{L^\infty(\Omega_n)}\lVert \nabla_\A \sigma_n \rVert_{L^2}^2 \right)\\
                                          &\leq c\lVert \nabla_\A \sigma_n \rVert_{L^2}^2,
  \end{align*}
  since $\sigma_n$  weakly converges in $W^{1,2}$, $\lVert \sigma_n \rVert_{W^{1,2}}$ is uniformly bounded, thus $\lVert \eta_n \rVert_{W^{1,2}}\leq c$.

  \step\label{step:An-convergence}Now, we prove that $A_n\to A_\infty$ strongly in $W^{1,2}(\Omega^1(\ad P_{G/\mathbb{Z}_2}))$. Note that, by $d \mathcal{H}(A_n,\sigma_n)(a_n-a_\infty)\to0$,
  \begin{align*}
    \lVert F_{A_n}-F_{A_\infty} \rVert_{L^2}^2&=\int_M \left\langle d(a_n-a_\infty),d(a_n-a_\infty) \right\rangle\\
                                              &=\int_M \left\langle da_n,d(a_n-a_\infty) \right\rangle+\left\langle da_\infty,d(a_n-a_\infty) \right\rangle\\
                                              &=\int_M\left\langle F_{A_n},d(a_n-a_\infty) \right\rangle+o(1)\\
                                              &=-\int_M\left\langle \nabla _{A_n}\sigma_n,(a_n-a_\infty)\sigma_n \right\rangle+o(1),
  \end{align*}
  Note also that on $M\setminus\Omega_n$, $\lvert \sigma_n \rvert \leq c_0$. We will continue to estimate 
  \begin{align*}
    \lVert F_{A_n}-F_{A_\infty} \rVert_{L^2}^2&\leq -\left( \int_{M\setminus \Omega_n}+\int_{\Omega_n} \right)\left( \left\langle \nabla _{A_n}\sigma_n,(a_n-a_\infty)\sigma_n \right\rangle \right)+o(1)\\
                                              &\leq  c_0 \int_{M\setminus\Omega_n}\left( \lvert \nabla _\A\sigma_n \rvert+\lvert a_n \rvert \right)\lvert a_n-a_\infty \rvert\\
                                              &\qquad+\int_{\Omega_n}\lvert \left\langle \nabla _{A_n}\sigma_n,\nu_n \right\rangle \rvert\lvert \sigma_n \rvert\lvert a_n-a_\infty \rvert+o(1)\\
                                              &\leq c\left( \lVert \nabla_\A \sigma_n \rVert_{L^2}+\lVert a_n \rVert_{L^2} \right)\lVert a_n-a_\infty \rVert_{L^2}\\
                                              &\qquad+c \lVert \left\langle \nabla _{A_n}\sigma_n,\nu_n \right\rangle \rVert_{L^2(\Omega_n)}\lVert \sigma_n \rVert_{L^4}\left( \lVert a_n\rVert_{L^4}  +\lVert a_\infty \rVert_{L^4} \right)+o(1),
  \end{align*}
  which tends to $0$ as $n\to\infty$. This is because, $a_n \rightharpoonup a_\infty$ in $L^4$ so $\lVert a_n \rVert_{L^4}$ is uniformly bounded and similarly $\lVert \sigma_n \rVert_{L^4}$, $\lVert \nabla_\A \sigma_n \rVert_{L^2}$ and $\lVert a_n \rVert_{L^2}$ are all uniformly bounded; and $a_n\to a_\infty$ strongly in $L^2$, $\left\langle \nabla _{A_n}\sigma_n,\nu_n \right\rangle\to0$ strongly in $L^2(\Omega_n)$ by the Claim. 

  In conclusion, we prove that $F_{A_n}\to F_{A_\infty}$ strongly in $L^2$, i.e., $da_n\to da_\infty$ strongly in $L^2$, since $a_n\to a_\infty$ strongly in $L^2$, we conclude that $A_n\to A_\infty$ strongly in $W^{1,2}$.

  \step \label{step:sigman-convergence} Finally, we prove that $\sigma_n\to \sigma_\infty$ strongly in $W^{1,2}(\Omega^1(E))$. Note that 
  \begin{align*}
    \lVert \nabla_\A \sigma_n-\nabla_\A \sigma_\infty \rVert_{L^2}^2
      &=\int_M\left\langle \nabla_\A \sigma_n,\nabla_\A (\sigma_n-\sigma_\infty) \right\rangle-\left\langle \nabla_\A \sigma_\infty,\nabla_\A (\sigma_n-\sigma_\infty) \right\rangle\\
      &=\int_M \left\langle \nabla_\A \sigma_n,\nabla_\A (\sigma_n-\sigma_\infty) \right\rangle+o(1)\\
      &=\int_M\left\langle \nabla _{A_n}\sigma_n-a_n\sigma_n,\nabla_{A_n} (\sigma_n-\sigma_\infty)-a_n(\sigma_n-\sigma_\infty) \right\rangle+o(1)\\
      &=\int_M\left\langle \nabla _{A_n}\sigma_n,\nabla _{A_n}(\sigma_n-\sigma_\infty) \right\rangle-\int_M\left\langle a_n\sigma_n,\nabla_\A (\sigma_n-\sigma_\infty) \right\rangle+o(1)\\
      &\qquad-\int_M\left\langle \nabla_\A \sigma_n,a_n(\sigma_n-\sigma_\infty) \right\rangle-\int_M\left\langle a_n\sigma_n,a_n(\sigma_n-\sigma_\infty) \right\rangle+o(1)\\
      &=\mathrm{I}-\mathrm{II}-\mathrm{III}-\mathrm{IV}+o(1).
  \end{align*}

  We will estimate the above four terms one by one. Firstly,  it is clear that 
  \begin{align*}
    -\int_M\left\langle \frac{1}{2}\mathfrak{T}(\sigma_n),\sigma_n-\sigma_\infty \right\rangle 
      &=-\int_M\sum_{\mu<\nu} \sum_a  \langle T_{\mu \nu a}\sigma_n, \sigma_n \rangle  \langle T_{\mu \nu a}\sigma_n, \sigma_n-\sigma_\infty \rangle\\
      &=-\int_M\sum_{\mu<\nu} \sum_a  \langle T_{\mu \nu a}\sigma_n, \sigma_n-\sigma_\infty \rangle  \langle T_{\mu \nu a}\sigma_n, \sigma_n-\sigma_\infty \rangle\\
      &\qquad-\int_M\sum_{\mu<\nu} \sum_a  \langle T_{\mu \nu a}\sigma_n, \sigma_\infty \rangle  \langle T_{\mu \nu a}\sigma_n, \sigma_n-\sigma_\infty \rangle\\
      &\leq C\int_M \lvert \sigma_n \rvert^2\lvert \sigma_\infty \rvert\lvert \sigma_n-\sigma_\infty \rvert,
  \end{align*}
  the last inequality follows from the fact that the first term is non-positive.
  Therefore,
  \begin{align*}
    \mathrm{I}&=\int_M\left\langle \nabla _{A_n}\sigma_n, \nabla _{A_n}(\sigma_n-\sigma_\infty) \right\rangle\\
              &=d \mathcal{H}(A_n,\sigma_n)(\sigma_n-\sigma_\infty) - \int_M\left\langle \mathrm{Ric}(\sigma_n),\sigma_n-\sigma_\infty \right\rangle-\frac{1}{2}\int_M\left\langle \mathfrak{T}(\sigma_n),\sigma_n-\sigma_\infty \right\rangle\\
              &= d \mathcal{H}(A_n,\sigma_n)(\sigma_n-\sigma_\infty) - \int_M\left\langle \mathrm{Ric}(\sigma_n-\sigma_\infty),\sigma_n-\sigma_\infty \right\rangle-\frac{1}{2}\int_M\left\langle \mathfrak{T}(\sigma_n),\sigma_n-\sigma_\infty \right\rangle+o(1)\\
              &\leq \lVert d \mathcal{H}(A_n,\sigma_n) \rVert_{W^{-1,2}}\lVert \sigma_n-\sigma_\infty \rVert_{W^{1,2}}
              +c\lVert \sigma_n-\sigma_\infty \rVert_{L^2}^2+c\lVert \sigma_n \rVert_{L^4}^2\lVert \sigma_n-\sigma_\infty \rVert_{L^2}+o(1),
  \end{align*}
  which tends to $0$ as $n\to\infty$, since $d \mathcal{H}(A_n,\sigma_n)\rightharpoonup 0$ in $W^{-1,2}$ and $\sigma_n-\sigma_\infty$ is $W^{1,2}$ bounded; and the last two terms tend to zero since $\sigma_n\to\sigma_\infty$ strongly in $L^2$.

  For the other three terms, we can use the fact that $A_n\to A_\infty$ in $W^{1,2}(\Omega^1(\ad P_{G/\mathbb{Z}_2}))$ as $n\to\infty$, which is already proved in \cref{step:An-convergence}. For example,
  \begin{align*}
    \mathrm{II}&=\int_M\left\langle a_n\sigma_n,\nabla_\A (\sigma_n-\sigma_\infty) \right\rangle\\
               &=\int_M\left\langle (a_n-a_\infty)\sigma_n,\nabla_\A (\sigma_n-\sigma_\infty \right\rangle+\left\langle a_\infty\sigma_n,\nabla_\A (\sigma_n-\sigma_\infty \right\rangle\\
               &\leq \lVert a_n-a_\infty \rVert_{L^4}\lVert \sigma_n \rVert_{L^4}\lVert \nabla_\A (\sigma_n-\sigma_\infty) \rVert_{L^2}
               +\int_M\left\langle a_\infty\sigma_n,\nabla_\A (\sigma_n-\sigma_\infty) \right\rangle\\
               &\leq \int_M\left\langle a_\infty\sigma_n,\nabla_\A (\sigma_n-\sigma_\infty) \right\rangle+o(1),
  \end{align*}
  this is because $a_n\to a_\infty$ strongly in $L^4$ and $\lVert \nabla_\A (\sigma_n-\sigma_\infty) \rVert_{L^2}\leq \lVert \sigma_n-\sigma_\infty \rVert_{W^{1,2}}\leq \lVert \sigma_n \rVert_{W^{1,2}}+\lVert \sigma_{\infty} \rVert_{W^{1,2}}$ is bounded. Thus, since $\alpha_\infty,\sigma_\infty\in L^4$ and $\sigma_n \rightharpoonup \sigma_\infty$ in $W^{1,2}$, we know that 
  \begin{align*}
    \mathrm{II} &\leq \int_M\left\langle a_\infty\sigma_n,\nabla_\A (\sigma_n-\sigma_\infty) \right\rangle+o(1)\\
                &=\int_M\left\langle a_\infty(\sigma_n-\sigma_\infty),\nabla_\A (\sigma_n-\sigma_\infty) \right\rangle+\left\langle a_\infty\sigma_\infty,\nabla_\A (\sigma_n-\sigma_\infty) \right\rangle+o(1)\\
                &=\int_M\left\langle a_\infty(\sigma_n-\sigma_\infty),\nabla_\A (\sigma_n-\sigma_\infty) \right\rangle+o(1)\\
                &\leq \lVert a_\infty \rVert_{L^\infty}\lVert( \sigma_n-\sigma_\infty) \rVert_{L^2}\lVert \nabla_\A (\sigma_n-\sigma_\infty) \rVert_{L^2}+o(1),
  \end{align*}
  which tends to zero as $n\to\infty$, because $a_\infty$ is smooth and $\sigma_n\to\sigma_\infty$ strongly in $L^2$ and $\lVert \nabla_\A (\sigma_n-\sigma_\infty) \rVert_{L^2}$ is bounded. 

  Hence, we can conclude that $\sigma_n$ strongly converges to $\sigma_\infty$ in $W^{1,2}(\Omega^1(E))$ as $n\to\infty$. Combining this with \cref{step:An-convergence}, we have completed the proof of \cref{thm:main}.
\end{proof}


\vspace*{2em}
\end{document}